\documentclass[english,11pt]{article}
\usepackage{babel,amsmath,amssymb,amsthm}   
\usepackage{enumerate}
\usepackage{enumitem}
\usepackage{tikz-cd}
\usepackage[matrix,arrow,tips,curve]{xy}
\usepackage{listings}
\usepackage{mathtools}
\usepackage[titletoc]{appendix}
\usepackage{url}
\usepackage{hyperref}
\usepackage{authblk}
\newcommand{\pr}{\mathbb{P}}
\newcommand{\bs}{\mathfrak{b}}
\newcommand{\I}{\mathfrak{I}}
\newcommand{\M}{\mathcal{M}}
\newcommand{\R}{\mathbb{R}}
\newcommand{\Bl}{\operatorname{Bl}}

\newcommand{\NE}{\operatorname{NE}}
\newcommand{\ma}{\mathcal}
\newcommand{\ol}{\mathcal{O}}

\newcommand{\w}{\widetilde}

\newcommand{\Pic}{\operatorname{Pic}}
\newcommand{\Bs}{\operatorname{Bs}}
\newcommand{\Sing}{\operatorname{Sing}}

\newcommand{\dg}{\operatorname{deg}}
\newcommand{\Aut}{\operatorname{Aut}}
\newcommand{\Id}{\operatorname{Id}}
\newcommand{\Eff}{\operatorname{Eff}}
\newcommand{\Nef}{\operatorname{Nef}}
\newcommand{\Mov}{\operatorname{Mov}}
\newcommand{\Ker}{\operatorname{Ker}}
\newcommand{\im}{\operatorname{Im}}

\newtheorem{thm}{Theorem}[section]

\newtheorem{lem}[thm]{Lemma}

\newtheorem{prop}[thm]{Proposition}
\newtheorem{cor}[thm]{Corollary}
\newtheorem{nota}[thm]{Notation}
\newtheorem{rem}[thm]{Remark}

\title{\textbf{Anticanonical geometry of the blow-up of $\pr^4$ in $8$ points and its Fano model}}
\author{Zhixin Xie}
\date{}
\begin{document}
\maketitle
\begin{abstract}
    Building on the work of Casagrande-Codogni-Fanelli, we develop our study on the birational geometry of the Fano fourfold $Y=M_{S,-K_S}$ which is the moduli space of semi-stable rank-two torsion-free sheaves with $c_1=-K_S$ and $c_2=2$ on a polarised degree-one del Pezzo surface $(S,-K_S)$. Based on the relation between $Y$ and the blow-up of $\pr^4$ in $8$ points, we describe completely the base scheme of the anticanonical system $|{-}K_Y|$. We also prove that the Bertini involution $\iota_Y$ of $Y$, induced by the Bertini involution $\iota_S$ of $S$, preserves every member in $|{-}K_Y|$. In particular, we establish the relation between $\iota_Y$ and the anticanonical map of $Y$, and we describe the action of $\iota_Y$ by analogy with the action of $\iota_S$ on $S$.
\end{abstract}
\section{Introduction}
Fano manifolds are classified up to dimension three. There
are $10$ deformation families of $2$-dimensional Fano manifolds, and $105$
deformation families of $3$-dimensional Fano manifolds (classified by Mori-Mukai and Iskovskih, see \cite{MR463151,MR503430,MR1022045,MR2112566}). 

In dimension $4$, toric Fano manifolds have been classified by Batyrev \cite{MR1703904} and Sato \cite{MR1772804}, and Fano manifolds of index $r>2$ have been classified. Among Fano fourfolds, those of index one have special positions: K\"{u}chle constructed a number of examples with Picard number one, and explained some known results with lists of related problems (see \cite{MR1454560}).
To find and classify Fano fourfolds of index one, Coates, Corti and others have embarked
on a program using mirror symmetry (\cite{MR3469127}, and see a list of examples in \cite{MR4101417}), where heavy computer calculations are
involved. A complete classification might not be desirable, but it is interesting to exhibit some Fano fourfolds with special geometric properties, for example, those with Picard number $\rho$ close to the conjectural boundary $\rho\leq 18$, or those whose anticanonical system has non-empty base locus. In order to study Fano manifolds with large Picard number (see \cite{MR3014481}), Casagrande introduced the invariant called \textit{Lefschetz defect}, and developed fruitful results in this direction (\cite{MR3092275,MR4297039}).

Let $Y\coloneqq M_{S,-K_S}$ be the moduli spaces of semi-stable rank-two torsion-free sheaves with $c_1=-K_S$, $c_2=2$ on a polarised degree-one del Pezzo surface $(S,-K_S)$. The moduli spaces $Y=M_{S,-K_S}$ form a remarkable family of smooth Fano fourfolds with Picard number $9$. The study of this family is motivated by two issues. Firstly,
for Fano fourfolds with large Picard number (e.g. at least $7$), only few examples which are not products of del Pezzo surfaces are known. As pointed out in \cite[Sect.~1,B]{MR3942925}, the family of Fano fourfolds $Y$ is the only known example of Fano fourfolds with Picard number at least $9$, which is not a product of surfaces.
Secondly, it is delicate to find examples of Fano fourfolds whose anticanonical system has non-empty base locus, since most Fano fourfolds classified so far are toric, which implies that any ample line bundle on them is globally generated. Some examples are constructed in \cite[Chapter 6.3]{heube2016} as complete intersections of two hypersurfaces in weighted projective spaces; two families are identified in \cite{secci2021fano} as Fano fourfolds with Picard number $3$ and having some contraction onto a smooth Fano threefold. In \cite[Thm.~1.10]{MR3942925}, it is shown that the base locus of the anticanonical system $|{-}K_Y|$ has positive dimension. Therefore, the geometry of $Y$ is worth detailed understanding.

\medskip
\noindent
{\bf A. The anticanonical system of the Fano model $Y$.}
The birational geometry of $Y= M_{S,-K_S}$ is related to the birational geometry of the blow-up $X$ of $\pr^4$ at $8$ points. In \cite[Lem.~5.18]{MR3942925}, an explicit relation between $X$ and $Y$ is given: the Fano fourfold $Y$ is obtained from $X$ by flipping the strict transforms of the lines through all pairs of blown up points and of the quartic curves through $7$ blown up points in $\pr^4$. Thanks to this relation between $X$ and $Y$, it is shown in \cite[Lem.~7.5, Cor.~7.6]{MR3942925} that the base locus of $|{-}K_Y|$ contains the strict transform $R_Y$ of a smooth rational quintic curve passing through the $8$ blown up points in $\pr^4$, and that $|{-}2K_Y|$ is base-point-free. We complete the study of the anticanonical system and show more precisely that:
\begin{thm}\label{mainbasecheme}
For the Fano fourfold $Y\coloneqq M_{S,-K_S}$, the base scheme of $|{-}K_Y|$ is the reduced smooth curve $R_Y$.
\end{thm}
As a direct application, we obtain the smoothness of a general member in the anticanonical system.
\begin{cor}\label{maincorsmooth}
Let $D\in|{-}K_Y|$ be a general divisor. Then $D$ is smooth.
\end{cor}

\medskip

\noindent
{\bf B. The Bertini involution of the Fano model $Y$.}
Now we turn our attention to the automorphism group of $Y$. In \cite[Sect.~4]{MR3942925}, a group morphism $\rho$ between the Picard groups of the degree-one del Pezzo surface $S$ and of the moduli spaces $Y=M_{S,-K_S}$ is defined. This morphism $\rho$ induces an isomorphism between the automorphism groups of $S$ and of $Y$ (see \cite[6.15]{MR3942925}). In particular, there is an involution $\iota_Y$ of $Y$ which is induced by the Bertini involution $\iota_S$ of $S$. 

We mention here that another motivation behind the study of the Bertini involution $\iota_Y$ is the understanding of the corresponding birational involutions $\iota_X$ of $X$ and $\iota_{\pr^4}$ of $\pr^4$. These birational maps $\iota_X$ and $\iota_{\pr^4}$ are classically known, as they can be defined via the Cremona action of the Weyl group $W(E_8)$ on sets of $8$ points in $\pr^4$ (see \cite{MR1007155} and \cite{MR661778}). Nevertheless, the classical definitions of $\iota_X$ and $\iota_{\pr^4}$ do not give a geometric description of these maps. In \cite[Prop.~8.9, Cor.~8.10]{MR3942925}, a factorisation of these maps is given as smooth blow-ups and blow-downs using the interpretation of $X$ as a moduli space of vector bundles on $S$. In view of the relation among $Y$, $X$ and $\pr^4$, understanding one of the involutions helps describe the behaviour of the others.

By the analogy of $Y$ and $S$, one expects that the action of $\iota_Y$ on $Y$ has similar properties as the action of $\iota_S$ on $S$, where the latter is well understood (see for example \cite[8.8.2]{MR2964027}). To emphasize their analogy, we recall that the Bertini involution $\iota_S$ on $S$ can be described as follows. The bianticanonical system $|{-}2K_S|$ is base-point-free and defines a 2:1-cover with image a quadric cone in $\pr^3$. The Bertini involution $\iota_S$ is then defined to be the associated covering involution.
By construction, the Bertini involution $\iota_S$ on $S$ preserves every element of $|{-}2K_S|$. Since a divisor $D\in |{-}K_S|$ defines an element $2D\in |{-}2K_S|$, we see that $\iota_S$ preserves every divisor in $|{-}K_S|$. In view of the abstract construction of $\iota_Y$ on the Fano fourfold $Y$, the same methode cannot be applied to decide whether $\iota_Y$ preserves every divisor in $|{-}K_Y|$. However, by analysing the anticanonical map of $Y$, we show that the same property holds for $Y$.
\begin{thm}\label{maininvantican}
The Bertini involution $\iota_Y$ of the Fano fourfold $Y\coloneqq M_{S,-K_S}$ preserves every divisor in $|{-}K_Y|$.
\end{thm}

To understand the Bertini involution $\iota_Y$ on the Fano fourfold, our approach is analysing its behaviour on a special surface $W_Y$ which is invariant by $\iota_Y$. This surface $W_Y$ is the strict transform of the cubic scroll swept out by the pencil of elliptic normal quintics in $\pr^4$ through the $8$ blown up points; in particular, it contains the curve $R_Y$. Inspired by the similarity with degree-one del Pezzo surfaces, we study the morphism defined by the restricted bianticanonical system of $Y$ on $W_Y$, and we give the following description of $\iota_Y$ restricted to $W_Y$.
\begin{prop}\label{mainintroBertini}
The Bertini involution $\iota_Y$ of the Fano fourfold $Y\coloneqq M_{S,-K_S}$ preserves the surface $W_Y$, and its restriction $\iota_Y|_{W_Y}$ on $W_Y$ is the biregular involution defined by the double covering \[\phi_{|{-}2K_Y|_{W_Y}}\colon W_Y \to V_{2,4}\subset\pr^7,\] where $V_{2,4}\simeq \mathbb{F}_2$ is a rational normal scroll of bidegree $(2,4)$. In particular, the Bertini involution $\iota_Y$ acts as the identity on the curve $R_Y$ and $\iota_Y$ induces an involution on each elliptic fibre $F_Y$ of $W_Y\to\pr^1$. 

Furthermore, there exists a smooth curve $R'\sim 3(R_Y+F_Y)$ of genus $4$ on the surface $W_Y$, such that $R'$ is disjoint from $R_Y$ and contained in the fixed locus of $\iota_Y$.
\end{prop}
Since $R_Y$ is contained in the fixed locus of the Bertini involution $\iota_Y$, the involution can be lifted to the blow-up $\tilde{Y}$ of $Y$ along the curve $R_Y$.
We establish the relation between the resolved anticanonical map and the lifted involution on $\tilde{Y}$ as follows.
\begin{thm}\label{maininvfact}
Let $\mu\colon\tilde{Y}\to Y$ be the blow-up of $Y\coloneqq M_{S,-K_S}$ along the base curve $R_Y$ of $|{-}K_Y|$, and $E$ be the exceptional divisor. Let $f\colon\tilde{Y}\to\pr(H^0(Y,\ol_Y(-K_Y))^{\vee})\simeq\pr^5$ be the morphism defined by the base-point-free linear system $|\mu^*(-K_Y)-E|$. Then $f$ has generically degree $4$ with image $Q$ a smooth quadric hypersurface in $\pr^5$, and $f$ contracts the strict transform of the surface $W_Y$ to a conic in $\pr^5$. Moreover, $f|_E\colon E\to f(E)$ is a finite birational morphism such that the image $f(E)$ has degree $4$ in $\pr^5$.

Furthermore, the Bertini involution $\iota_Y$ of $Y$ can be lifted to $\tilde{Y}$, and the lifted involution $\iota_{\tilde{Y}}$ acts as the identity on $E$. Moreover, $f$ factors through the quotient $\tilde{Y}/\iota_{\tilde{Y}}$: 
\begin{center}
\begin{tikzcd}
\tilde{Y} \arrow[rd] \arrow[d, "\mu"'] \arrow[dd, "f"', bend right] &                                        \\
Y \arrow[d, dashed]                                                 & \tilde{Y}/\iota_{\tilde{Y}} \arrow[ld] \\
Q\subset\mathbb{P}^5                                                &                                       
\end{tikzcd}
\end{center}
\end{thm}

As open questions, one may like to understand the quotient $\tilde{Y}/\iota_{\tilde{Y}}$, the geometric interpretation of $\tilde{Y}/\iota_{\tilde{Y}}\to Q$, and to describe completely the fixed locus of $\iota_{\tilde{Y}}$ (see Remark \ref{fixedlocustildeY}). 

\medskip

\noindent
{\bf Plan.} We briefly explain the organisation of the paper. In Section \ref{prelimFano}, we summarise some results in \cite{MR3942925}, including the geometry of the Fano model $Y\coloneqq M_{S,-K_S}$, the connection between the blow-up $X$ of $\pr^4$ at $8$ points and the degree-one del Pezzo surface $S$, and the relation between $X$ and $Y$. We finish by recalling some basic properties of the Bertini involution of a degree-one del Pezzo surface.

In Section \ref{sec:antican}, we investigate the anticanonical system $|{-}K_Y|$ and the bianticanonical system $|{-}2K_Y|$. We prove Proposition \ref{mainbasecheme} by an additional analysis on the simplicial facets of the cone of effective divisors on $Y$. We also give some auxiliary results on $|{-}K_Y|$ and $|{-}2K_Y|$, which serve as key ingredients in the study of the Bertini involution of $Y$.

In Section \ref{Bert}, we study the action of the Bertini involution of $Y$. Subsection \ref{sec:bertini} is devoted to the proof of Proposition \ref{mainintroBertini}. We study the morphism defined by the bianticanonical system $|{-}2K_Y|$ restricted to the surface $W_Y$. Computations by {\tt Macaulay2} show that the image of $W_Y$ is a surface of degree $6$ in $\pr^7$, which helps us to describe completely the morphism; in particular, the morphism is finite of degree $2$ and gives an involution on the surface $W_Y$. By examining the action of this covering involution, we show that it coincides with the Bertini involution $\iota_Y$ restricted to the surface $W_Y$.

In Subsection \ref{sec:bertiniantican}, we study the geometry of the anticanonical map of $Y$.
Computations by {\tt Macaulay2} show that the image of $Y$ by the antincanonical map is a smooth quadric hypersurface $Q$ in $\pr^5$. We are then ready to prove Theorem \ref{maininvantican}.
The strategy is to prove by contradiction: we suppose that $\iota_Y$ does not preserve every divisor in $|{-}K_Y|$. We show that in this case, $\iota_Y$ induces a non-trivial involution $\iota_Q$ on $Q$. We then obtain a contradiction by analysing the fixed locus of the induced involution $\iota_Q$ and by studying the geometry of a special sub-linear systems of $|{-}K_Y|$ consisting of divisors containing the surface $W_Y$. Theorem \ref{maininvfact} is obtained as a consequence of the study to prove Theorem \ref{maininvantican}.

In Appendix \ref{appen}, we include the code for several computations in Section \ref{sec:antican} and Section \ref{Bert} using the software system {\tt Macaulay2}.

\medskip
{\bf Acknowledgements.} This project was initiated during my
stay in Turin. I heartily thank Cinzia Casagrande for her hospitality and fruitful conversations.
I would like to express my sincere gratitude to my supervisor, Andreas H\"{o}ring, for his patient guidance, his constant support and valuable suggestions. I also thank Daniele Faenzi for his help on computations, and Susanna Zimmermann for interesting discussions. 

Many results in this paper are based on computations using {\tt Macaulay2}. I would like to thank the developers for making their software open-source.
I thank the IDEX UCA JEDI project (ANR-15-IDEX-01) and the MathIT project for providing financial support.
\section{Preliminaries}\label{prelimFano}
We fix $S$ a general del Pezzo surface of degree one. Let $M_{S,L}$ be the moduli space of semi-stable (with respect to $L\in \Pic(S)$ ample) rank-two torsion free sheaves $\mathcal{F}$ on $S$ with $c_1(\mathcal{F})=-K_S$ and $c_2(\mathcal{F})=2$. Then it follows from the classical properties of the determinant line bundle that for the polarisation $L=-K_S$, the moduli space $Y\coloneqq M_{S,-K_S}$ is Fano.

For the degree-one del Pezzo surface $S$, we introduce the following notions (see \cite[Sect.~2.1]{MR3942925}). A conic on $S$ is a smooth rational curve such that $-K_S\cdot C=2$ and $C^2=0$. Every such conic yields a conic bundle $S\to\pr^1$ having $C$ as fibre. There are $2160$ conics (as classes of a curve) in $H^2(S,\mathbb{Z})$. A big divisor $h$ on $S$ which realises $S$ as the blow-up $\sigma:S\to\pr^2$ at $8$ distinct points is called a cubic. We have $h=\sigma^*\ol_{\pr^2}(1)$. There are $17280$ cubics (as classes of a curve) in $H^2(S,\mathbb{Z})$.

\begin{nota}\label{notadP}
Given a cubic $h$, we use the following notation:
\begin{itemize}
\item
 $\sigma_h\colon S\to \pr^2$ is the birational map defined by $h$ 
\item $q_1,\dotsc,q_8\in\pr^2$ are the points blown up by $\sigma_h$
\item $e_i\subset S$ is the exceptional curve over $q_i$, for $i=1\dotsc,8$
\item $C_i\subset S$ is the  transform of a general line through $q_i$, so that $C_i\sim h-e_i$, for $i=1,\dotsc,8$
\item $\ell_{ij}\subset S$ is the transform of the line $\overline{q_iq_j}\subset\pr^2$, so that $\ell_{ij}\sim h-e_i-e_j$ and $\ell_{ij}$ is a $(-1)$-curve, for $1\leq i<j\leq 8$.
\end{itemize}
\end{nota}

\subsection{The Fano model $Y$}
By \cite[Prop.~1.6]{MR3942925}, the moduli space $Y\coloneqq M_{S,-K_S}$ is a smooth, rational Fano fourfold with index one and Picard number $9$.
For such a moduli space $Y$, the determinant map $\rho\colon H^2(S,\R)\to H^2(Y,\R)$ is an isomorphism (see \cite[Thm.~1.3]{MR3942925}) and yields a completely explicit description of the relevant cones of divisors $\Eff(Y)$, $\Mov(Y)$ and $\Nef(Y)$, as well as the cone of effective curves $\NE(Y)$. We cite the following statements for the cone of effective divisors $\Eff(Y)$ and the cone of effective curves $\NE(Y)$, and refer the readers to \cite[Sect.~6]{MR3942925} for the description of the other relevant cones.
\begin{prop}[\cite{MR3942925}, Sect.~2.3, Cor.~6.2]\label{Eff}
The determinant map $\rho\colon H^2(S,\R)\to H^2(Y,\R)$ yields an isomorphism between $\mathcal{E}$ and $\Eff(Y)$, where $\mathcal{E}$ is the subcone of $\Nef(S)$ generated by the conics:
\[
\mathcal{E}\coloneqq\langle C\,|\,C\text{ a conic}\rangle\subset H^2(S,\R).
\]
Hence, the cone $\Eff(Y)$ has $2160$ extremal rays, each generated by a fixed divisor $E_C$, where $C\subset S$ is a conic.

Moreover, given a cubic $h$, $(2h+K_S)^{\perp}\cap \mathcal{E}$ is a simplicial facet (i.e. a face of codimension one) of $\mathcal{E}$, generated by the conics $C_i$ for $i=1,\dots,8$ (notations as in Notation \ref{notadP}). Hence, the fixed divisors $E_{C_i}$ for $i=1,\dots,8$ generate a simplicial facet of $\Eff(Y)$.
\end{prop}

\begin{prop}[\cite{MR3942925}, Prop.~1.7]\label{smallloci}
The cone of effective curves $\NE(Y)$  
has $240$ extremal rays, and is isomorphic to  $\NE(S)$. If $\ell$ is a $(-1)$-curve on $S$, the corresponding extremal ray of $\NE(Y)$ is generated by the class of a line $\Gamma_{\ell}$ in $P_{\ell}\cong\pr^2\subset Y$.
The corresponding elementary contraction is a small contraction, sending $P_{\ell}$  to a point. \end{prop}

The determinant map $\rho$ also relates the two automorphism groups $\Aut(Y)$ and $\Aut(S)$.
By \cite[Thm.~1.9]{MR3942925},
the map $\psi\colon\Aut(S)\to\Aut(Y)$ given by $\psi(\phi)[\mathcal{F}]=[(\phi^{-1})^*\mathcal{F}]$, for $\phi\in\Aut(S)$ and $[\mathcal{F}]\in Y$, is a group isomorphism. In particular, $\Aut(Y)$ is finite; if $S$ is general,
then \[\Aut(Y)=\{\Id_Y,\iota_Y\},\] where  $\iota_Y\colon Y\to Y$ is induced by the Bertini involution of $S$. 
We still call the involution $\iota_Y\coloneqq\psi(\iota_S)$ of $Y$ the Bertini involution. Explicitly, $\iota_Y\colon Y\to Y$ is given (see \cite[Def.6.19]{MR3942925}) by
$\iota_Y([\mathcal{F}])=[\iota_S^*\mathcal{F}]$. We have a commutative diagram:
\begin{equation}\label{bertini}\begin{gathered}{\footnotesize
\xymatrix{{H^2(S,\R)}\ar[r]^{\iota_S^*}\ar[d]_{\rho}&{H^2(S,\R)}\ar[d]^{\rho}\\{H^2(Y,\R)}\ar[r]^{\iota_Y^*}&{H^2(Y,\R).}
 }}\end{gathered}\end{equation}

Finally, motivated by the analogy with del Pezzo surface of degree one, the study of the base loci of the anticanonical and the bianticanonical linear systems of $Y$ gives the following:
\begin{thm}[\cite{MR3942925}, Thm.~1.10]\label{system}
The linear system $|{-}K_Y|$ has a base locus of positive dimension, while
the linear system $|{-}2K_Y|$ is base point free.
\end{thm}

\subsection{The blow-up $X$ of $\pr^4$ at $8$ general points}
\subsubsection{Degree one del Pezzo surfaces and blow-ups of $\pr^4$ in $8$ points}
For $S=\Bl_{q_1,\dots,q_8}\pr^2$ and $X=\Bl_{p_1,\dots,p_8}\pr^4$ the blow-ups respectively of $\pr^2$ and $\pr^4$ at $8$ general points, there is a classical association between these two varieties due to Gale duality. The following is summarised from \cite[2.21]{MR3942925}; for further details of the association, we refer to \cite[2.18]{MR3942925}.

Let $h$ be a cubic on $S$. We associate to $(S,h)$ a blow-up $X$ of $\pr^4$ in $8$ points in  general linear position as follows.

Let $q_1,\dotsc,q_8\in\pr^2$ be the points blown up under the birational morphism 
$S\to\pr^2$ defined by $h$ (the points $q_1,\dots,q_8$ are in general linear position by \cite[Rem.~2.20]{MR3942925}), and let $p_1,\dotsc,p_8\in\pr^4$ be the associated points to $q_1,\dotsc,q_8\in\pr^2$ (the points $p_1,\dots,p_8$ are in  general linear position by \cite[Lem.~2.19]{MR3942925}). Then we set \[X=X_h=X_{(S,h)}\coloneqq\Bl_{p_1,\dotsc,p_8}\pr^4.\]
We always assume that $q_1,\dotsc,q_8\in\pr^2$ and  $p_1,\dotsc,p_8\in\pr^4$ are associated as ordered sets of point.

Conversely, let $X$ be a blow-up of $\pr^4$ in $8$ general points. Differently from the case of surfaces, the blow-up map $X\to\pr^4$ is unique and thus $X$ determines $p_1,\dotsc,p_8\in\pr^4$ up to projective equivalence. The $8$ points $p_1,\dots,p_8\in\pr^4$ in turn determine 
$q_1,\dotsc,q_8\in\pr^2$ up to projective equivalence, and thus a pair $(S,h)$ such that $X\cong X_{(S,h)}$. The pair $(S,h)$ is unique up to isomorphism, therefore $S$ is determined up to isomorphism, and $h$ is determined up to
the action of the automorphism group $\Aut(S)$ on cubics.

\subsubsection{Notation for the blow-up $X$ of $\pr^4$ at $8$ points}\label{notationP^4}
Let $p_1,\dotsc,p_8\in\pr^4$ be $8$ points in general linear position, and set $X\coloneqq\Bl_{p_1,\dotsc,p_8}\pr^4$. We use the following notation:
\begin{itemize}
\item $E_i\subset X$ is the exceptional divisor over $p_i\in\pr^4$, for $i=1,\dotsc,8$
\item $H\in\Pic(X)$ is the pull-back of $\ol_{\pr^4}(1)$
\item $L_{ij}\subset X$ is the transform of the line $\overline{p_ip_j}\subset\pr^4$,  for $1\leq i<j\leq 8$
\item $e_i\subset E_i$ is a line, for $i=1,\dotsc,8$
\item $\gamma_i\subset\pr^4$ is the rational normal quartic  through $p_1,\dotsc,\check{p}_i,\dotsc,p_8$, for $i=1,\dotsc,8$
\item $\Gamma_i\subset X$ is the transform of  $\gamma_i\subset\pr^4$, for $i=1,\dotsc,8$.
\end{itemize}

\subsection{From the blow-up $X$ to the Fano model $Y$}
We recall the explicit relation between $X$ and $Y$:
\begin{lem}[\cite{MR3942925}, Lem.~5.18]\label{sequence}
The birational map $\xi\colon X\dasharrow Y$ is the composition of $36$ ($K$-positive) flips: first the flips of  $L_{ij}$ for $1\leq i<j\leq 8$, and then the flips of   $\Gamma_k$ for $k=1,\dotsc,8$. There is a commutative diagram:
{\footnotesize$$\xymatrix{&{\widehat{X}}\ar[dl]\ar[dr]&\\
{X}\ar@{-->}[rr]^{\xi}&&Y}$$}where $\widehat{X}\to X$ is the blow-up of the curves $L_{ij}$  and $\Gamma_k$, with every exceptional divisor isomorphic to $\pr^1\times\pr^2$ with normal bundle $\ma{O}(-1,-1)$, and $\widehat{X}\to Y$ is the blow-up of $36$ pairwise disjoint smooth rational surfaces.
\end{lem}
\begin{nota}
We use the following notation:
\begin{itemize}
    \item $P_{\ell_{ij}}\subset Y$ is the flipped surface replacing $L_{ij}\subset X$, for $1\leq i<j\leq 8$
    \item $P_{e_k}\subset Y$ is the flipped surface replacing $\Gamma_k\subset X$, for $k=1,\dots,8$.
\end{itemize}
\end{nota}

We will sometimes write $\xi_h\colon X_h\dasharrow Y$ to stress that $X_h$ and $\xi_h$ depend on the chosen cubic $h$ (while $Y$ does not).
Denote by $\eta_h$ the composition map:

\begin{center}
\begin{tikzcd}
Y \arrow[r, "\xi_h^{-1}"', dashed] \arrow[rr, "\eta_h", dashed, bend left] & X_h \arrow[r] & \mathbb{P}^4
\end{tikzcd}
\end{center}

\subsection{The Bertini involution of $S$}
We recall some basic properties of the Bertini involution of a del Pezzo surface of degree one.

\begin{prop}[\cite{MR2964027},Thm.~8.3.2]
Suppose that $S$ is a del Pezzo surface of degree $1$. Then
\begin{enumerate}[label=\normalfont(\roman*)]
    \item $|{-}K_S|$ is a pencil of genus $1$ curves with smooth general member and one base point;
    \item $|{-}2K_S|$ is base-point-free and defines a morphism $\phi_{|{-}2K_S|}\colon S\to\pr^3$ which is finite of degree $2$ with image $Q$ a quadric cone.
\end{enumerate}
The Bertini involution $\iota_S\colon S\to S$ is the biregular involution defined by the double covering \[\phi_{|{-}2K_S|}\colon S\to Q.\]
\end{prop}
For $S$ general, $\iota_S$ is the unique non-trivial automorphism of $S$.
The pull-back $\iota_S^*$ acts on $\Pic(S)$ (and on $H^2(S,\R)$) by  fixing
 $K_S$ and acting as $-1$ on $K_S^{\perp}$ (see \cite[\S 8.8.2]{MR2964027}). This yields:
\begin{equation}
\label{w_0}
\iota_S^*\gamma=2(\gamma\cdot K_S)K_S-\gamma\quad\text{ for every $\gamma\in H^2(S,\R)$.}
\end{equation}
The fixed locus of $\iota_S$ is a smooth irreducible curve of genus $4$ isomorphic to the branch curve of the double cover and the base point of $|{-}K_S|$. The fixed curve belongs to the linear system $|-3K_S|$.

\section{Anticanonical and bianticanonical linear systems of the Fano model $Y$}\label{sec:antican}

Let $S$ be a degree-one del Pezzo surface, and $Y\coloneqq M_{S,-K_S}$ be the associated Fano fourfold.
To analyse the anticanonical linear system $|{-}K_Y|$, we introduce a special surface as follows.
\begin{lem}[\cite{MR3942925}, Lem.~7.2]\label{scroll}
Let $p_1,\dotsc,p_8\in\pr^4$ be general points. 
Then there is a pencil  of elliptic normal quintics in $\pr^4$ through $p_1,\dotsc,p_8$, which sweeps out a cubic scroll $W\subset\pr^4$.

Let moreover 
$q_1,\dotsc,q_8\in\pr^2$ be the associated points to  $p_1,\dotsc,p_8\in\pr^4$. Then there is a birational map $\alpha\colon W\to\pr^2$ such that  $\alpha(p_i)=q_i$ for $i=1,\dotsc,8$, $\alpha$ sends the pencil  of elliptic normal quintics to the pencil of plane cubics through $q_1,\dotsc,q_8$, and $\alpha$ is the blow-up of the ninth base point $q_0\in\pr^2$ of the pencil of plane cubics.
\end{lem}
Let $W'\subset X$ be the transform of the cubic scroll $W\subset\pr^4$. By \cite[(7.3)]{MR3942925}, we have the following diagram:
\begin{equation}\label{W}\begin{gathered}{\footnotesize
\xymatrix{&{W'\subset X}\ar[dl]_{\eta}\ar[dr]\ar[dd]^{\alpha'}&\\
{W\subset\pr^4}\ar[dr]_{\alpha}&&S\ar[dl]^{\sigma}\\
&{\pr^2}&
}}\end{gathered}\end{equation}
where $\eta\colon W'\to W$ is the blow-up of $p_1,\dotsc,p_8$, so the composition $\alpha'\coloneqq\alpha\circ
\eta\colon W'\to\pr^2$ is the blow-up of $q_0,\dotsc,q_8$. Thus
$W'$ is isomorphic to the blow-up of $S$ in the base point of  $|{-}K_S|$. Hence, there is an elliptic fibration $\pi\colon W'\to\pr^1$,  where the smooth fibres are the transforms of the elliptic normal quintics through
$p_1,\dotsc,p_8$ in $\pr^4$, and every fibre is integral.

\begin{lem}[\cite{MR3942925}, Lem.~7.4]\label{locus}
The surface $W'\subset X$ is disjoint from $L_{ij}$ for $1\leq i<j\leq 8$ and from $\Gamma_k$ for $k=1,\dotsc,8$, and $W'$ is contained in the open subset where $\xi\colon X\dasharrow Y$ is an isomorphism. 
\end{lem}
We denote by $W_Y$ the strict transform of $W'$ in $Y$. Then $W_Y\simeq W'$.

\begin{lem}[\cite{MR3942925}, Lem.~7.5, Lem.~7.7, Rem.~7.10]\label{fibration}
We have $(-K_X)_{|W'}=\ol_{W'}(R+2F)$ and $R=\Bs|(-K_X)_{|W'}|$,
where $F\subset W'$ is a fibre of the elliptic fibration, and $R\subset W'$ is a $(-1)$-curve and a section of the elliptic fibration. The curves $R$ and $F$ satisfy $-K_X\cdot R=-K_X\cdot F=1$ and $E_i\cdot R=E_i\cdot F=1$ for every $i=1,\dots,8$, so $R\equiv F$ in $X$ and $\xi(R)\equiv\xi(F)$ in $Y$.

Moreover, let $R_4\subset \pr^4$ be the images of $R$ under $\eta\colon W'\subset X\to W\subset \pr^4$ (see diagram \eqref{W}).
Then $R_4$ is a smooth rational quintic curve through $p_1,\dotsc,p_8$
\end{lem}

\begin{cor}[\cite{MR3942925}, Cor.~7.6]\label{conclusion}
The base locus of $|{-}K_X|$ contains the smooth rational curve $R$, and the base locus of $|{-}K_Y|$ contains the
 smooth rational curve $\xi(R)$. 
\end{cor}
We denote by $R_Y$ the smooth rational curve $\xi(R)$ contained in the base locus of $|{-}K_Y|$, and $F_Y$ a fibre of the elliptic fibration $W_Y\to\pr^1$.

\begin{lem}\label{normalbundle}
The normal bundle $\ma{N}_{R_Y/Y}\cong\ol_{\pr^1}(-1)\oplus\ol_{\pr^1}^{\oplus 2}.$
\end{lem}
\begin{proof}
Since $R_4$ is a rational quintic curve in $\pr^4$, one has \[\ma{N}_{R_4/\pr^4}\cong\ol_{\pr^1}(a)\oplus\ol_{\pr^1}(b)\oplus\ol_{\pr^1}(c)\] with $a\leq b\leq c$ and $a+b+c=23$.
Since $\mathcal{T}_{\pr^4}|_{R_4}\twoheadrightarrow \mathcal{N}_{R_4/\pr^4}\to 0$, one has that $\mathcal{N}_{R_4/\pr^4}$ is ample. Hence, we deduce that $a,b,c>0$.
Moreover, by {\tt Macaulay2} (see Listing \ref{lst:normal}), \[ h^0(R_4,\ma{N}^*_{R_4/\pr^4}\otimes \mathcal{O}_{\pr^4}(1)\otimes\omega^*_{R_4})=1,\] 
we deduce that $a=7$ and $b,c>7$. Hence, $b=c=8$.
Therefore, by \cite[B.6.10]{MR1644323}, one has
\[
\ma{N}_{R/X}\cong\ol_{\pr^1}(-1)\oplus\ol_{\pr^1}^{\oplus 2}.
\]
As $R$ is disjoint from the indeterminacy locus of the map $\xi_h$, we deduce that
\[
\ma{N}_{R_Y/Y}\cong\ol_{\pr^1}(-1)\oplus\ol_{\pr^1}^{\oplus 2}.
\]
\end{proof}

\begin{rem}(see also \cite[Remark 7.8]{MR3942925})\label{baseP4}
In $\pr^4$, let $\M$ be the linear system of quintic hypersurfaces with multiplicity at least $3$ at $8$ general points. Then by {\tt Macaulay2} (see Listing \ref{lst:base}) the base ideal $\bs(\M)$ is the intersection of the ideals of $28$ line $\overline{p_ip_j}$ for $1\leq i<j\leq 8$, the ideals of $8$ rational normal quartic curves $\gamma_k$ for $k=1,\dots,8$ and the ideal of the rational quintic curve $R_4$. This shows that the base scheme of $|{-}K_Y|$, in the open subset of $Y$ where $\eta_h:Y\dashrightarrow\pr^4$ is an isomorphism, is $R_Y$ (reduced) minus the $8$ points of intersection with the $8$ exceptional divisors $\xi(E_1),\dots,\xi(E_8)$. The base locus of $|{-}K_Y|$ is thus given by $R_Y$, possibly union some other components contained in $\xi(E_1),\dots,\xi(E_8)$.
\end{rem}

\begin{lem}\label{basedisjoint}
The base locus of the anticanonical system $|{-}K_Y|$ is disjoint from the surfaces $P_{\ell_{ij}}$ and $P_{e_k}$, for $1\leq i<j\leq 8$ and $k=1,\dots,8$.
\end{lem}
\begin{proof}
Consider the commutative diagram in Lemma \ref{sequence}:
\[
\begin{tikzcd}
                            & \hat{X} \arrow[ld, "p"'] \arrow[rd, "q"] &   \\
X \arrow[rr, "\xi", dashed] &                                          & Y
\end{tikzcd}
\]
where $p\colon\hat{X}\to X$ is the blow-up of $X$ along the curves $L_{ij}$ and $\Gamma_k$ with every exceptional divisor isomorphic to $\pr^1\times\pr^2$, and $q\colon\hat{X}\to Y$ is the blow-up of $36$ pairwise disjoint smooth rational surfaces $P_{\ell_{ij}}$ and $P_{e_k}$, for $1\leq i<j\leq 8$ and $k=1,\dots,8$.

Suppose by contradiction that there exists a base point $y$ of $|{-}K_Y|$ contained in some flipped surface that we denote by $P$ (which is one of the surfaces $P_{\ell_{ij}}$ or $P_{e_k}$). Denote by $C\subset X$ the corresponding flipping curve (which is one of the curves $L_{ij}$ or $\Gamma_k$). 

Let $E$ be the sum of exceptional divisors over $L_{ij}$ for $1\leq i<j\leq 8$ and over $\Gamma_k$ for $k=1,...,8$.
Since
\[
p^*(-K_X)-2E = -K_{\hat{X}}= q^*(-K_Y)-E,
\]
one has
\[
q^*(-K_Y)=p^*(-K_X)-E.
\]
Let $E_y\simeq\pr^1$ be the exceptional fibre in $\hat{X}$ above $y$. Then $E_y$ is contained in $\Bs |q^*(-K_Y)|=\Bs |p^*(-K_X)-E|$ and $E_y$ is mapped surjectively onto $C$.

Since the blow-up of the $8$ points $X=\Bl_{p_1,\dots,p_8}\to\pr^4$ is an isomorphism near a general point of $C$, the base scheme of $|{-}K_X|$ is generically reduced along $C$ by Remark \ref{baseP4}. Hence, the linear system $|p^*(-K_X)-E|$ is base-point-free above the generic point of $C$. This contradicts the fact that $\Bs |p^*(-K_X)-E|$ contains a curve which is mapped surjectively onto $C$.
\end{proof}

\begin{rem}
More generally, the proof of Lemma \ref{basedisjoint} shows the following. Let $X,Y$ be smooth projective fourfolds. Let $\xi\colon X\dashrightarrow Y$ be an anti-flip. In \cite[Thm.~1.1]{MR1006374}, Kawamata showed that for smooth projective fourfolds, there exists only one type of flip and it is obtained by blowing up a $\pr^2$ with normal bundle $\ol_{\pr^2}(-1)^{\oplus 2}$ (the exceptional locus of the blowing up is $\pr^2\times\pr^1$) and
blowing down this $\pr^2\times\pr^1$ to $\pr^1$. Thus $\xi$ (anti-)flips a smooth curve $C\subset X$ to a smooth surface $P\subset Y$. If $\Bs|{-}K_X|$ is reduced in the generic point of $C$, then $|{-}K_Y|$ is base-point-free on $P$.
\end{rem}

\begin{cor}\label{basecurve}
The curve $R_Y$ is the unique base curve in $\Bs|{-}K_Y|$ of anticanonical degree $1$. Therefore, $R_Y$ is independent of the choice of cubic $h$.
\end{cor}
\begin{proof}
Let $C\subset \Bs|{-}K_Y|$ be a base curve contained in some exceptional divisor $\xi(E_i)$, for $i=1,\dots,8.$ Let $\tilde{C}$ be its strict transform in $X$. By Lemma \ref{basedisjoint}, the curve $C$ is disjoint from the indeterminacy locus of $\xi^{-1}$. Hence, one has $-K_Y\cdot C=-K_X\cdot\tilde{C}$ and $\tilde{C}\subset E_i.$

Since $-K_X=5H-3\sum_{j=1}^8 E_j$, $H\cdot\tilde{C}=0,$ $E_j\cdot\tilde{C}=0$ for $j\neq i$, and $E_i\cdot\tilde{C}\leq -1$, one has \[-K_Y\cdot C=-K_X\cdot\tilde{C}\geq 3.\]
Therefore, the curve $R_Y$ is the unique base curve satisfying $-K_Y\cdot R_Y=-K_X\cdot R=1$.
\end{proof}

\begin{cor}\label{uniquefixed}
Let $B\subset Y$ be an irreducible component of the (set-theoretic) base locus of $|{-}K_Y|$, which is distinct from $R_Y$. Then for every simplicial facet $\langle E_{C_1},\dots, E_{C_8} \rangle$ of $\Eff(Y)$ (notation as in Notation \ref{notadP} and Proposition \ref{Eff}), there exists a unique $E_{C_i}$ for $i=1,\dots,8$ such that $B\subset E_i$.
\end{cor}
\begin{proof}
Given a cubic $h$, consider the simplicial facet $\langle E_{C_1},\dots, E_{C_8} \rangle$ of $\Eff(Y)$, where $C_i\sim h-e_i$ for $i=1,\dots,8$ (notation as in Notation \ref{notadP}). Then $E_{C_i}$ are the strict transforms of the exceptional divisors $E_i\simeq \pr^3\subset X_h= X$ under $\xi_h = \xi\colon X\dashrightarrow Y$.

Since $B$ is distinct from $R_Y$, we deduce that $B$ is contained in some fixed divisor $E_{C_i}$ by Remark \ref{baseP4}. By the construction of the composition of flips $\xi$ (see Lemma \ref{sequence}), the intersection of two fixed divisors $E_{C_j}$ and $E_{C_k}$ (for $k\neq j$) is the union of the flipped surfaces $P_{\ell_{jk}}$ and $P_{e_l}$ for $l\neq j,k$. Hence, by Lemma \ref{basedisjoint}, the fixed divisor $E_i$ containing $B$ is unique.
\end{proof}

\begin{proof}[Proof of Theorem \ref{mainbasecheme}]
We first show that $R_Y$ is the unique irreducible component of the (set-theoretic) base locus of $|{-}K_Y|$.

Let $h$ be a cubic. Let $C_i$ be a conic such that $C_i\sim h-e_i$ for $i=1,\dots,8$ (notation as in Notation \ref{notadP}). Let $E_i\coloneqq E_{C_i}$, where we use the notation of Proposition \ref{Eff}. By the same proposition, $E_1,\dots,E_8$ generate a simplicial facet of $\Eff(Y)$. Suppose by contradiction that there exists another component $B$ distinct from $R_Y$ of the base locus of $|{-}K_Y|$. Then by Corollary \ref{uniquefixed}, we may suppose that $B\subset E_1$ and $B\not\subset E_2, E_3,\dots,E_8.$

Let $i,j,k,l$ be distinct indices in $\{ 1,\dots, 8\}$. Consider the conics $C'_l$ such that $C'_l\sim 2h-e_i-e_j-e_k-e_l$ and the corresponding fixed divisors $F_{ijkl}\coloneqq E_{2h-e_i-e_j-e_k-e_l}$.

{\em Claim}. The fixed divisors $E_i,E_j,E_k$ and $F_{ijkl}$ for $l\in\{1,\dots,8\}$ distinct from $i,j,k$ generate a simplicial facet of $\Eff(Y)$.

Indeed, by Proposition \ref{Eff}, it is enough to find a cubic $h'$ such that $2h'+K_S$ is orthogonal to the $8$ conics $C_i,C_j,C_k$ and $C'_l$ for $l\in\{1,\dots,8\}$ distinct from $i,j,k$.

We take $h'\sim 2h-e_i-e_j-e_k$. Then we can check that
\begin{align*}
    C_i &\sim h'-\ell_{jk}\\
    C_j &\sim h'-\ell_{ik}\\
    C_k &\sim h'-\ell_{ij}\\
    C'_l&\sim h'-e_l
\end{align*}
and $2h'+K_S$ is orthogonal to the above $8$ conics. Moreover, $h'$ is nef and big, and the corresponding birational map $\sigma_{h'}\colon S\to\pr^2$ contracts the $8$ pairwise disjoint $(-1)$-curves $\ell_{jk},\ell_{ik},\ell_{ij},e_l$ for $l\neq i,j,k$. Hence, $h'$ is a cubic. This proves the claim.

\medskip

We will repeatedly use Corollary \ref{uniquefixed} in the following.
\begin{itemize}
\item Consider the simplicial facet generated by $E_1,E_2,E_3,F_{1234},F_{1235},F_{1236},F_{1237},F_{1238}$. Then $B\not\subset F_{1234},F_{1235},F_{1236},F_{1237},F_{1238}.$
\item Consider the simplicial facet generated by $E_2,E_3,E_4,F_{1234},F_{2345},F_{2346},F_{2347},F_{2348}$. Then $B$ is contained in one of the fixed divisors $F_{2345},F_{2346},F_{2347},F_{2348}$. We may suppose that $B\subset F_{2345}.$ Then $B\not\subset F_{2346},F_{2347},F_{2348}.$
\item Consider the simplicial facet generated by $E_2,E_3,E_5,F_{1235},F_{2345},F_{2356},F_{2357},F_{2358}$. Then $B\not\subset F_{2356},F_{2357},F_{2358}$.
\item Consider the simplicial facet generated by $E_2,E_3,E_6,F_{1236},F_{2346},F_{2356},F_{2367},F_{2368}$. Then by what precedes, we know that $B$ is contained in one of the fixed divisors $F_{2367},F_{2368}$. We may suppose that $B\subset F_{2367}$. Then $B\not\subset F_{2368}$.
\item Consider the simplicial facet generated by $E_2,E_3,E_7,F_{1237},F_{2347},F_{2357},F_{2367},F_{2378}$. Then $B\not\subset F_{2378}$.
\item Finally, consider the simplicial facet generated by $E_2,E_3,E_8,F_{1238},F_{2348},F_{2358},F_{2368},F_{2378}$. Then by what precedes, we know that $B$ is contained in none of these $8$ fixed divisors, which contradicts Corollary \ref{uniquefixed}.
\end{itemize}
Therefore, the curve $R_Y$ is the unique irreducible component of the base locus of $|{-}K_Y|$.

Now we show that the base scheme of $|{-}K_Y|$ is the reduced curve $R_Y$, i.e. there are no embedded points. Indeed, given a cubic $h$, consider the birational map $\eta_h\colon Y\dashrightarrow\pr^4$. By Remark \ref{baseP4}, the base scheme of $|{-}K_Y|$ is the reduced curve $R_Y$ with possible embedded points which have support in the $8$ points of intersection with the $8$ exceptional divisors of $\eta_h$. By varying $h$, we may consider another map $\eta_{h'}\colon Y\dashrightarrow\pr^4$ with other $8$ exceptional divisors, so that we get $8$ different points of intersection on $R_Y$. Such a cubic $h'$ exists because otherwise, there is a base point $y$ on $R_Y$ such that for every simplicial facet $\langle E_{C_1},\dots, E_{C_8} \rangle$ of $\Eff(Y)$ the point $y$ is contained in a unique $E_{C_i}$, and thus we obtain a contradiction by replacing $B$ with $y$ in the above paragraph. Hence, there is no embedded base point on $R_Y$.
\end{proof}

\begin{proof}[Proof of Corollary \ref{maincorsmooth}]
Since the base scheme $\Bs |{-}K_Y|$ is the smooth curve $R_Y$ by Proposition \ref{mainbasecheme}, we can apply \cite[Prop.~6.8]{MR1102273} which implies that a general member in $|{-}K_Y|$ is smooth.
\end{proof}

In the rest of this section, we collect some auxiliary results which will be used in Section \ref{Bert}.
\begin{lem}\label{sec}
For a general point $x\in R_4$ (notation as in Lemma \ref{fibration}), there exists a unique divisor in $\mathcal{M}$ which has multiplicity $3$ at $x$: it is the secant variety of the elliptic normal quintic through the nine points $p_1,\dots,p_8$ and $x$.

By varying $x$ on $R_4$, one obtains a one-dimensional family $Sec$ of divisors in $\mathcal{M}$ with scheme-theoretic intersection $\Bs Sec$ defined by the ideal $\mathfrak{b}(Sec)$.  Then
\[
(\bs(Sec):\bs(\M)):\I_W=\I_W,
\]
where the scheme defined by the ideal $\I_W$ is the reduced surface $W$.
\end{lem}
\begin{proof}
We choose a random point $x$ on $R_4$ which is not one of the $8$ blown up points. Let $\M_{x,3}$ be the linear subspace of divisors in $\M$ having multiplicity at least $3$ at the point $x$. Then $\dim \M_{x,3}=0$ by ${\tt Macaulay2}$ (see Listing \ref{lst:secant}) and thus the unique element in $\M_{x,3}$ is the secant variety $Sec(E_x)$, where $E_x$ is the elliptic normal quintic in $W$ passing through the point $x$ and the $8$ blown up points.

Let $Sec$ be the family of secant varieties $Sec(E_x)$ for $x$ varying on $R_4$ and $\bs(Sec)$ be the ideal associated to the scheme-theoretic intersection $\Bs Sec$ of the family $Sec$. Let $\bs_3(Sec)$ be the ideal associated to the scheme-theoretic intersection of three general secant varieties in $Sec$ (obtained by choosing three distinct random points on $R_4$ and intersecting the corresponding secant varieties).

By ${\tt Macaulay2}$ (see Listing \ref{lst:W}), the quotient $\I_S\coloneqq (\bs_3(Sec):\bs (\M))$ has degree $6$ and dimension $2$. Let $\I_W$ be the ideal of singular locus of the variety defined by $\I_S$. Then by ${\tt Macaulay2}$ (see Listing \ref{lst:W}), $\I_W$ has dimension $2$ and degree $3$; moreover, the variety defined by $\I_W$ is smooth and one has $(\I_S:\I_W) = \I_W.$ Since each of these secant varieties in $Sec$ contains the cubic scroll $W$, we deduce that the variety defined by $\I_W$ is indeed the surface $W$.

Let $\M_W$ be the sub-linear system of effective divisors in $\M$ containing the surface $W$. By {\tt Macaulay2} (see Listing \ref{lst:anticanW}), the base ideal $\bs(\M_W)$ is equal to $\bs_3(Sec)$.
Since $Sec$ is a family of divisors in $\mathcal{M}_W$, we deduce that $\bs_3(Sec)=\bs(Sec)=\bs(\M_W)$.
\end{proof}

\begin{lem}\label{Pldisj}
The surface $W_Y$ is unique, i.e. $W_Y$ is independent of the choice of cubic $h$. Therefore, $W_Y$ is disjoint from every one of the loci $P_\ell$ of the small extremal rays of $Y$.
\end{lem}
\begin{proof}
Let $Sec_Y$ be the family of the strict transforms in $Y$ of the secant varieties in $Sec$. Let $\M_{Y,3}$ be the family of divisors in $|{-}K_Y|$ having multiplicity $3$ at some point on $R_Y$. Then the two families $\M_{Y,3}$ and $Sec_Y$ are equal, as $\dim \M_{x,3}=0$ for a general point $x\in R_4$ by Lemma \ref{sec} and $\eta_h$ is an isomorphism at the generic point of $R_4$.

Suppose by contradiction that $W_Y$ depends on $h$. Then there exist two distinct surfaces $W_{Y,h}$ and $W_{Y,h'}$. Let $\Bs \M_{Y,3}$ be the scheme-theoretic intersection of the family $\M_{Y,3}$. Then by Lemma \ref{sec}, one has the following set-theoretic inclusion:
\[
\Bs M_{Y,3} \supset W_{Y,h}\cup W_{Y,h'}.
\]
Since $W_{Y,h'}$ contains the curve $R_Y$ which is generically in the locus where $\xi^{-1}_h\colon Y\dashrightarrow X_h$ is an isomorphism, we deduce that $W_{Y,h'}$ is not contracted by $\xi^{-1}_h$.

Since the surface $\xi^{-1}_h(W_{Y,h'})$ contains the curve $R$, this surface cannot be contained in any exceptional locus $E_i$, $i=1,...,8$ of $X_h\to \pr^4$, and thus it cannot be contracted; we denote by $W_{h'}$ its image in $\pr^4$.
Therefore, $\Bs Sec$ contains two distinct surfaces $W$ and $W_{h'}$, which contradicts Lemma \ref{sec}.

Since by Lemma \ref{locus} the surface $W_Y$ is disjoint from the indeterminacy locus of the map $\xi^{-1}_h\colon Y\dasharrow X_h$, which is a union of some of the loci $P_\ell$ (depending on $h$), and $W_Y$ is the same for all $h$, we deduce that $W_Y$ is disjoint from every one of the loci $P_\ell$.
\end{proof}

\begin{lem}\label{restriction}
\begin{enumerate}[label=\normalfont(\roman*)]
\item We have $h^0(W_Y,\ol_{W_Y}(-K_Y))=3$. The restriction \[r_1\colon H^0(Y,\ol_Y(-K_Y))\to H^0(W_Y,\ol_{W_Y}(-K_Y))\] is surjective.
\item We have $h^0(W_Y,\ol_{W_Y}(-2K_Y))=8$. The restriction \[r_2\colon H^0(Y,\ol_Y(-2K_Y))\to H^0(W_Y,\ol_{W_Y}(-2K_Y))\] is surjective.
\end{enumerate}
\end{lem}
\begin{proof}
Since $-K_{W_Y}\sim F_Y$ and $-K_Y|_{W_Y}\sim R_Y+2F_Y$ by Lemma \ref{fibration}, by the Riemann-Roch formula one has $\chi(W_Y,-K_Y|_{W_Y})=3$. Since $-K_Y|_{W_Y}$ is ample on $W_Y$ and $-K_{W_Y}$ is nef, by Kodaira vanishing theorem one has \[h^j(W_Y,\ol_{W_Y}(-K_Y))=h^j(W_Y,\ol_{W_Y}(K_{W_Y}-K_{W_Y}+(-K_Y)))=0\] for $j=1,2$. Therefore, $h^0(W_Y,\ol_{W_Y}(-K_Y))=3$. 
The same argument can be applied to obtain $h^0(W_Y,\ol_{W_Y}(-2K_Y))=8$. 

(i) By {\tt Macaulay2} (see Listing \ref{lst:anticanW}), 
\[
h^0(\pr^4,\ol_{\pr^4}(5)\otimes\mathcal{I}_{p_1,\dots,p_8}^3\otimes\mathcal{I}_{W})=3.
\]
Since $H^0(Y,\ol_{Y}(-K_Y))\simeq H^0(\pr^4,\ol_{\pr^4}(5)\otimes\mathcal{I}_{p_1,\dots,p_8}^3)$, and the surface $W_Y$ is disjoint from the indeterminacy locus of $\eta_h$ by Lemma \ref{locus} and $W_Y$ is not contained in the exceptional locus of $\eta_h$, we deduce that
\[
H^0(Y,\ol_{Y}(-K_Y)\otimes\mathcal{I}_{W_Y})\simeq H^0(\pr^4,\ol_{\pr^4}(5)\otimes\mathcal{I}_{p_1,\dots,p_8}^3\otimes\mathcal{I}_{W}).
\]
Hence,
\[ h^0(Y, \ol_{Y}(-K_Y)\otimes\mathcal{I}_{W_Y})=3. \] 
As $h^0(Y,\ol_Y(-K_Y))=6$ and $h^0(W_Y,\ol_{W_Y}(-K_Y))=3$, we deduce that the restriction morphism
\[
H^0(Y,\ol_Y(-K_Y))\to H^0(W_Y,\ol_{W_Y}(-K_Y))
\]
is surjective.

(ii) By {\tt Macaulay2} (see Listing \ref{lst:bianticanW}), 
\[
h^0(\pr^4,\ol_{\pr^4}(10)\otimes\mathcal{I}_{p_1,\dots,p_8}^6\otimes\mathcal{I}_{W})=21.
\]
Since $H^0(Y,\ol_{Y}(-2K_Y))\simeq H^0(\pr^4,\ol_{\pr^4}(10)\otimes\mathcal{I}_{p_1,\dots,p_8}^6)$ and by the same argument as above, we deduce that
\[
H^0(Y,\ol_{Y}(-2K_Y)\otimes\mathcal{I}_{W_Y})\simeq H^0(\pr^4,\ol_{\pr^4}(10)\otimes\mathcal{I}_{p_1,\dots,p_8}^6\otimes\mathcal{I}_{W}).
\]
Hence,
\[ h^0(Y, \ol_{Y}(-2K_Y)\otimes\mathcal{I}_{W_Y})=21. \] 
As $h^0(Y,\ol_Y(-2K_Y))=29$ and $h^0(W_Y,\ol_{W_Y}(-2K_Y))=8$, we deduce that the restriction morphism
\[
H^0(Y,\ol_Y(-2K_Y))\to H^0(W_Y,\ol_{W_Y}(-2K_Y))
\]
is surjective.
\end{proof}

\section{The Bertini involution of the Fano model $Y$}\label{Bert}
Let $S$ be a degree-one del Pezzo surface, and $Y\coloneqq M_{S,-K_S}$ be the associated Fano fourfold.
In this section, we study the action of the Bertini involution $\iota_Y$ on the Fano fourfold $Y$, which is analogous to the action of the Bertini involution $\iota_S$ on the surface $S$. We first notice that by the diagram (\ref{bertini}) and the behaviour of $\iota_S$ described in (\ref{w_0}), the invariant part of $H^2(Y,\R)$ by the action of $\iota_Y$ is $\R K_Y$.
\subsection{Action of the Bertini involution on the surface $W_Y$}\label{sec:bertini}
In this subsection, we further our study of the involution $\iota_Y$ by looking at its action on the surface $W_Y$ (which is the strict transform of the cubic scroll swept out by the pencil of elliptic normal quintics in $\pr^4$). The aim of this subsection is to prove Proposition \ref{mainintroBertini}.

We start by showing that the surface $W_Y$ is invariant by the Bertini involution $\iota_Y$.
\begin{lem}\label{invY}
The Bertini involution $\iota_Y$ preserves the curve $R_Y$ and the surface $W_Y$. Moreover, $(\iota_Y|_{W_Y})^*(e_i)\sim -2K_Y|_{W_Y}-e_i$ and $(\iota_Y|_{W_Y})^*(F_Y)\sim F_Y$, where $e_i$ is the exceptional curve of $\eta_h|_{W_Y}\colon W_Y\to W$ for $i=1,\dots,8$.
\end{lem}
\begin{proof}
Since $\iota_Y$ preserves the family of divisors in the anticanonical system $|{-}K_Y|$, the involution $\iota_Y$ preserves the base locus of $|{-}K_Y|$. Thus $\iota_Y(R_Y)=R_Y$ by Proposition \ref{mainbasecheme}.

Let $x$ be a general point in $R_Y$. Then by Lemma \ref{sec}, there exists a unique divisor in $|{-}K_Y|$ having multiplicity $3$ at $x$: it is the strict transform in $Y$ of the secant variety of the elliptic normal quintic through $p_1,\dots,p_8$ and $\eta_h(x)$ in $\pr^4$. In particular, this divisor has multiplicity $3$ along the elliptic fibre of $W_Y$ through $x$. By varying $x$ in $R_Y$, this gives a one-dimensional family $\M_{Y,3}$ of divisors in $|{-}K_Y|$, which is preserved by $\iota_Y$. On the other hand, the intersection of these divisors is the surface $W_Y$, so $W_Y$ is preserved by $\iota_Y$. Let $D_1\in M_{Y,3}$ and $D_2=\iota_Y(D_1)\in M_{Y,3}$. Let $F_1$ (resp. $F_2$) be the elliptic fibre of $W_Y$ along which $D_1$ (resp. $D_2$) has multiplicity $3$. Then $\iota_Y(F_1)=F_2$, and thus $\iota_Y$ preserves the family of elliptic fibres of $W_Y$, i.e. $(\iota_Y|_{W_Y})^*(F_Y)\sim F_Y$.

By \cite[7.12]{MR3942925}, one has $\iota_Y^*(\xi(E_i))\sim -2K_Y-\xi(E_i)$. Hence, $(\iota_Y|_{W_Y})^*(e_i)\sim -2K_Y|_{W_Y}-e_i.$
\end{proof}

Now we investigate the morphism defined by the linear system $|{-}2K_Y|_{W_Y}|$.
\begin{prop}\label{invW}
The linear system $|{-}2K_Y|_{W_Y}|$ defines a finite morphism $\phi\colon W_Y\to V\subset\pr^7$ of degree $2$, where $V=V_{2,4}\simeq\mathbb{F}_2$ is a rational normal scroll of bidegree $(2,4)$. There is a non-trivial involution $i$ of $W_Y$ such that $\phi=\phi\circ i$. Moreover, $i$ is the identity on $R_Y$ and $i$ induces an involution on each elliptic fibre of $W_Y$.
\end{prop}

\begin{proof}
Since $h^0(W_Y,\ol_{W_Y}(-2K_Y))=8$ (see Lemma \ref{restriction}), and $|{-}2K_Y|$ is base-point-free by Theorem \ref{system}, the linear system $|{-}2K_Y|_{W_Y}|$ defines a morphism $\phi\colon W_Y\to V\subset\pr^7$, where $V$ is the image of $W_Y$.

\medskip
\noindent
{\em Claim}. $V$ is a surface of degree $6$ in $\pr^7$, the image of an elliptic fibre $F_Y$ by $\phi$ is a line and the image of $R_Y$ by $\phi$ is a conic.

Since the restriction morphism $H^0(Y,\ol_Y(-2K_Y))\to H^0(W_Y,\ol_{W_Y}(-2K_Y))$ is surjective by Lemma \ref{restriction} $(ii)$, the restriction of the morphism $\phi_{|{-}2K_Y|}$ defined by $|{-}2K_Y|$ to the surface $W_Y$ coincides with the morphism $\phi$, i.e. $\phi=\phi_{|{-}2K_Y|}|_{W_Y}$.

In $\pr^4$, let $2\M$ be the linear system of hypersurfaces of degree $10$ with multiplicity at least $6$ at the $8$ general points $p_1,...,p_8$. Consider the map $\phi_{2\M}$ defined by the linear system $2\M$. Then by {\tt Macaulay2} (see Listing \ref{lst:image}), the image of the surface $W$ by $\phi_{2\M}$ is a surface of degree $6$, the image of an elliptic normal quintic through the $8$ points by $\phi_{2\M}$ is a line and the image of the rational quintic $R_4$ through the $8$ points by $\phi_{2\M}$ is a conic.
This proves the claim.

\medskip
Since $(-2K_Y|_{W_Y})^2=4(R_Y+2F_Y)^2=12$, and the image of $W_Y$ by $\phi$ is of degree $6$, we deduce that $\phi$ is of degree $2$. As $-K_Y$ is ample, the morphism $\phi$ does not contract any curve and thus it is a finite morphism of degree $2$.

Since the linear system $|{-}2K_Y|_{W_Y}|$ has no fixed divisor, the image $V$ is not contained in any hyperplane of $\pr^7$ (see for example \cite[II.6]{MR1406314}), i.e. $V$ is non-degenerate. Hence, $V$ is a non-degenerate irreducible surface of degree $6$ (variety of minimal degree) in $\pr^7$, and by \cite[p.~525]{MR1288523} we deduce that $V$ is a rational normal scroll $V_{k,l}$ of bidegree $(k,l)$, with $0\leq k\leq l$ and $k+l=6$. In particular, $V$ is isomorphic to one of the following: a cone over a rational normal curve of degree $6$, $\pr^1\times\pr^1$, or a Hirzebruch surface $\mathbb{F}_{l-k}$, where the minimal section is mapped to the rational normal curve of degree $k$, and the fibres are mapped to lines. Therefore, $\phi$ is a finite morphism between two normal surfaces and by \cite[(2.3)]{MR695657}, there is a non-trivial involution $i$ of $W_Y$ such that $\phi=\phi\circ i$ and $V\simeq W_Y/i$.

Since the restriction of $\phi$ to a general fibre $F_Y$ induces a finite morphism from an elliptic curve to a line $l\subset V$, which cannot be an isomorphism, we deduce that $\phi^{-1}(l)=F_Y$ as $\phi$ is of degree $2$. Hence, $i$ induces an involution on $F_Y$.

Since $-2K_Y|_{W_Y}$ is $i$-invariant, one has $2(i^*(-K_Y|_{W_Y})-(-K_Y|_{W_Y}))\sim 0$. As $\Pic(W_Y)$ is torsion-free (this is because $W_Y$ is isomorphic to $\pr^2$ blown up at $9$ points), we deduce that $i^*(-K_Y|_{W_Y})\sim -K_Y|_{W_Y}$.
Since $R_Y$ is the base locus of $|{-}K_Y|_{W_Y}|$, the curve $R_Y$ is preserved by $i$. We claim that $R_Y$ is contained in the ramification locus of $\phi$. Indeed, suppose that $R_Y$ is not contained in the ramification locus of $\phi$. Then there exists a curve $C\subset V$ such that $R_Y=\phi^*(C)$. As $R_Y$ is a $(-1)$-curve on $W_Y$, one has
\[
-1=R_Y^2=(\phi^*(C))^2=\dg\phi\cdot C^2,
\]
i.e. $C^2=-\frac{1}{2}$. Hence, $C$ is not Cartier on $V$, i.e. $V$ is singular. In view of the classification of minimal degree varieties, we see that $V$ is a cone. But there is no curve with negative self-intersection number on a cone, which leads to a contradiction. Therefore, $R_Y$ is in the ramification locus. As $\phi$ is a double cover, we deduce that $i$ is the identity on $R_Y$.

Let $C=\phi(R_Y)$. Since $R_Y$ is contained in the ramification locus of $\phi$, and every point in $R_Y$ has ramification index $2$, one has
\[
R_Y^2=(\frac{1}{2}\phi^*(C))^2=\frac{1}{2}C^2.
\]
Since $R_Y$ is a $(-1)$-curve on $W_Y$, one has $C^2=-2$. Therefore, $V=V_{2,4}\simeq\mathbb{F}_2$, and $\phi(R_Y)$ is minimal section of $\mathbb{F}_2$ which is a conic.
\end{proof}

\begin{rem}\label{remfixedlocus}
Since $\phi$ is a finite morphism of degree $2$ between smooth surfaces, the ramification locus is a smooth divisor on $W_Y$ (see \cite[(2.5)]{MR695657}). Let $e$ be the minimal section of $V\simeq\mathbb{F}_2$ and $f$ be a fibre of $V$. Let $D$ be the ramification divisor. Then
\[
K_{W_Y}\sim\phi^*(K_S)+D.
\]
As $K_{W_Y}\sim -F_Y=-\phi^*(f)$, and $K_S\sim-2e-4f$, one has
\[
D\sim\phi^*(2e+3f).
\]
Let $B\subset V$ be the branch locus. Then $D=\frac{1}{2}\phi^*B$ and thus $B\sim 4e+6f$. As $e$ is contained in the branch locus, we can write $B=e+B_1$, where $B_1$ is a smooth curve disjoint from $e$. Then $B_1\sim 3e+6f$. Notice that $B_1$ is irreducible. Indeed, suppose that $B_1$ has at least two disjoint irreducible components. Then we can decompose $B_1$ as
\[
B_1 \sim (e + bf) + (2e + (6-b)f)
\]
with $0\leq b \leq 6$ and $(e + bf)\cdot(2e + (6-b)f)=0$. Hence $b=-2$, which leads to a contradiction.

Hence $D=R_Y + R'$, where $R'\sim\frac{1}{2}\phi^*(3e+6f)=3(R_Y+F_Y)$ is a smooth curve of genus $4$ which is disjoint from $R_Y$.
\end{rem}

Finally, we compare the action of the two automorphisms $i$ and $\iota_Y|_{W_Y}$ on $W_Y$.
\begin{lem}\label{invWe}
Let $e_i$ be the exceptional curves of $\eta_h|_{W_Y}\colon W_Y\to W$ for $i=1,\dots,8$. Then
\[
i^*(e_i)\sim -2K_Y|_{W_Y}-e_i.
\]
\end{lem}
\begin{proof}
For $i=1,\dots,8$, by {\tt Macaulay2} (see Listing \ref{lst:special}), there exists a unique hypersurface of degree $10$ with multiplicity at least $7$ at the point $p_i$ and multiplicity at least $6$ at $p_j$ for $j\neq i$. Moreover, this hypersurface does not contain the surface $W$. Therefore, the linear system $|{-}2K_Y|_{W_Y}-e_i|$ is non-empty.

Let $R_i\in |{-}2K_Y|_{W_Y}-e_i|$. Since $-K_Y|_{W_Y}\sim R_Y+2F_Y$, and $R_Y\cdot e_i=F_Y\cdot e_i=1$, one has $R_i^2=-1$, and  $R_i\cdot F_Y=R_i\cdot R_Y=1$. Hence, $R_i$ is a $(-1)$-curve on $W_Y$.

Since $e_i+R_i\in|{-}2K_Y|_{W_Y}|=\phi^*|\ol_V(1)|$, one has $R_i\sim i^*(e_i)$.
\end{proof}

\begin{prop}\label{coincide}
The involution $i$ coincides with the restriction of the Bertini involution $\iota_Y$ on the surface $W_Y$, i.e. $\iota_Y|_{W_Y}=i.$
\end{prop}
\begin{proof}
We first show that $(\iota_Y|_{W_Y})^*=i^*$. By Lemma \ref{invW}, Lemma \ref{invWe} and Lemma \ref{invY}, it is enough to show that $R_Y, F_Y$ and $e_i$ for $i=1,\dots, 8$ form a basis of $H^2(W_Y,\mathbb{R})$.

Since $W'$ is disjoint from the indeterminacy locus of $\xi_h$, it is equivalent to show that $R,F$ and $e_i$ for $i=1,\dots, 8$ form a basis of $H^2(W',\mathbb{R})$. We have the following diagram (see \eqref{W}):

\begin{center}
\begin{tikzcd}
                                           & W'\subset X \arrow[ld, "\eta"'] \arrow[dd, "\alpha'"] \\
W\subset\mathbb{P}^4 \arrow[rd, "\alpha"'] &                                                       \\
                                           & \mathbb{P}^2                                         
\end{tikzcd}
\end{center}
where $\alpha$ is the blow-up of $\pr^2$ at one point and $\eta$ is the blow-up of $W$ at $p_1,\dots,p_8$. Moreover, let $e_0\subset W$ be the $(-1)$-curve and $f_0\subset W$ be a fibre of the $\pr^1$-bundle on $W$, then by Lemma \ref{scroll} and Lemma \ref{fibration}, one has $F\sim \eta^*(2e_0+3f_0)-\sum_{i=1}^8 e_i$ and $R\sim \eta^*(e_0+4f_0)-\sum_{i=1}^8 e_i$. Therefore, $R,F$ and $e_i$ for $i=1,\dots,8$ form a basis of $H^2(W',\mathbb{R})$.

We have a group homomorphism $\rho_1\colon\Aut(W_Y)\to \Aut(H^2(W_Y,\R))$ given by $g\mapsto (g^{-1})^*$. Let $\Aut (R_Y,W_Y)$ be the subgroup of automorphisms in $\Aut(W_Y)$ fixing the curve $R_Y$. We show that the restriction $\rho_1\colon\Aut (R_Y,W_Y)\to \Aut(H^2(W_Y,\R))$ is injective, which implies that $\iota_Y|_{W_Y}=i$.

Since $R_Y$ is a $(-1)$-curve on $W_Y$, by blowing down $R_Y$, we obtain a rational surface $S'$ with $(-K_{S'})^2=1$, and the curve $R_Y$ is contracted to a point $x_0\in S'$. We denote by $\beta:W_Y\to S'$ the blow-up of $S'$ at $x_0$. Since $-K_{W_Y}$ is nef, we obtain that $-K_{S'}$ is nef by the projection formula (see for example \cite[Appendix A, A4]{MR0463157}). Moreover, since every fibre of $W_Y\to\pr^1$ is integral, there is no $K_{S'}$-trivial curve. Hence, $S'$ is a del Pezzo surface of degree one. By \cite[Prop.~8.2.39]{MR2964027}, the homomorphism $\rho_2\colon\Aut(S')\to \Aut H^2(S',\R)$ is injective. 

Let $\Aut(x_0,S')$ be the subgroup of automorphisms in $\Aut(S')$ fixing the point $x_0$. Then $\Aut(x_0,S')\simeq \Aut(R_Y,W_Y)$. Since $\Pic(W_Y)\simeq \beta^*\Pic(S')\oplus \mathbb{Z}[R_Y]$, the image $\rho_1(\Aut(R_Y,W_Y))$ is contained in a subgroup $G_1$ of $\Aut(H^2(W_Y,\R))$ such that $G_1\simeq\Aut(H^2(S',\R))$. Hence, we have the following diagram:

\begin{center}
\begin{tikzcd}
{\operatorname{Aut}(R_Y,W_Y)} \arrow[rr, "\rho_1"] \arrow[d, "\simeq"] &  & G_1 \arrow[d, "\simeq"]                  \\
{\operatorname{Aut}(x_0,S')} \arrow[rr, "\rho_2"]                      &  & {\operatorname{Aut}(H^2(S',\mathbb{R}))}
\end{tikzcd}  
\end{center}
Since $\rho_2$ is injective, the restriction $\rho_1\colon\Aut(R_Y,W_Y)\to G_1\subset \Aut(H^2(W_Y,\R))$ is injective.
\end{proof}

\begin{proof}[Proof of Proposition \ref{mainintroBertini}]
The first paragraph follows from Lemma \ref{invY}, Proposition \ref{invW} and Proposition \ref{coincide}. The second paragraph follows from Remark \ref{remfixedlocus}.
\end{proof}

\subsection{Action of the Bertini involution on the anticanonical system}\label{sec:bertiniantican}
In this subsection, we study the action of the involution $\iota_Y$ on the anticanonical system $|{-}K_Y|$. This is closely related to the anticanonical map of $Y= M_{S,-K_S}$.
\begin{lem}\label{blpR}
Let $\mu\colon\tilde{Y}\to Y$ be the blow-up of $Y$ along the curve $R_Y$ which is the base scheme of $|{-}K_Y|$. Let $E$ be the exceptional divisor and $\tilde{D}$ be the strict transform of a general member $D\in|{-}K_Y|$. Then $|\tilde{D}|=|\mu^*(-K_Y)-E|$ is base-point-free and induces a morphism $f\colon\tilde{Y}\to\pr(H^0(Y,\ol_Y(-K_Y))^{\vee})\simeq\pr^5$ with image $Q$ a smooth quadric hypersurface, and $f$ has generically degree $4$. We have the following commutative diagram:
\begin{equation}\label{resbase}\begin{gathered}{\footnotesize
\xymatrix{&{\tilde{Y}}\ar[dl]_{\mu}\ar[dr]^f&\\
{Y}\ar@{-->}[rr]^{\phi_{|{-}K_Y|}}&&Q\subset\pr^5}}\end{gathered}\end{equation}
The following statements hold:
\begin{enumerate}[label=\normalfont(\roman*)]
\item The Bertini involution $\iota_Y$ can be lifted to an involution $\iota_{\tilde{Y}}$ of $\tilde{Y}$, which preserves the exceptional divisor $E$ and induces an involution on each $\pr^2$ above a point of $R_Y$.

\item The Bertini involution $\iota_Y$ induces an involution $\iota_{\pr^5}$ of $\pr^5$, which preserves the quadric hypersurface $Q$. Denote by $\iota_Q$ its restriction on $Q$. Then $\iota_Q\circ f=f\circ\iota_{\tilde{Y}}$.
\end{enumerate}
\end{lem}
\begin{proof}
In $\pr^4$, let $\M$ be the linear system of quintic hypersurfaces with multiplicity at least $3$ at $8$ general points. Then by {\tt Macaulay2}, the image of $\pr^4$ by the map defined by $\M$ is a smooth quadric hypersurface $Q$ in $\pr^5$.

Let $E$ be the exceptional divisor of $f$. Since
    \[\mu^*(-K_Y)^4=(-K_Y)^4=13,\] \[\mu^*(-K_Y).E^3=-K_Y.R_Y=1,\] \[\mu^*(-K_Y)^3.E=\mu^*(-K_Y)^2.E^2=0,\] \[E^4=-K_Y.R_Y+2g(R_Y)-2 =-1,\]
    one has $\w{D}^4=8$. Hence $\phi_{|{-}K_Y|}$ (and also $f$) has generically degree $4$.
    
(i) Follows from the fact that $R_Y$ is contained in the fixed locus of $\iota_Y$ (see Proposition \ref{mainintroBertini}).

(ii) The pull-back $\iota_Y^*$ induces an involution on $H^0(-K_Y,\ol_Y(-K_Y))$, and thus an involution of $\pr(H^0(Y,\ol_Y(-K_Y))^{\vee})\simeq\pr^5$ preserving $\phi_{|{-}K_Y|}(Y)=Q$.

Let $s\in H^0(\tilde{Y},\ol_{\tilde{Y}}(\tilde{D}))$ be a global section which is zero at the point $\iota_{\tilde{Y}}(x)$, where $x$ is a point in $\tilde{Y}$. Then for $s'\coloneqq \iota_{\tilde{Y}}^*(s)\in H^0(\tilde{Y},\ol_{\tilde{Y}}(\iota_{\tilde{Y}}^*\tilde{D}))\simeq H^0(\tilde{Y},\ol_{\tilde{Y}}(\tilde{D}))$, one has
\[
s'(x)=(\iota_{\tilde{Y}}^*(s))(x)=s(\iota_{\tilde{Y}}(x))=0.
\]
Hence, \[\phi_{|\tilde{D}|}(\iota_{\tilde{Y}}(x))=\{ s\in H^0(\tilde{Y},\ol_{\tilde{Y}}(\tilde{D})) \mid s(\iota_{\tilde{Y}}(x))=0\},\]
\[\phi_{|\iota_{\tilde{Y}}^*\tilde{D}|}(x)=\{ s'\in H^0(\tilde{Y},\ol_{\tilde{Y}}(\iota_{\tilde{Y}}^*\tilde{D})) \mid s'(x)=0\}.\]
Therefore, we obtain the following commutative diagram:
\begin{center}
\begin{tikzcd}
{\tilde{Y}} \arrow[rr, "f"] \arrow[d, "\iota_{\tilde{Y}}"] &  & {\pr(H^0(\tilde{Y},\ol_{\tilde{Y}}(\iota_{\tilde{Y}}^*(\tilde{D})))^{\vee})} \arrow[d, "\iota_{\pr^5}"]                  \\
{\tilde{Y}} \arrow[rr, "f"]                      &  & {\pr(H^0(\tilde{Y},\ol_{\tilde{Y}}(\tilde{D}))^{\vee}).}
\end{tikzcd}  
\end{center}
Thus, $\iota_Q\circ f=f\circ\iota_{\tilde{Y}}$.
\end{proof}

\begin{rem}
The following statements are equivalent:
\begin{enumerate}[label=\normalfont(\alph*)]
\item The Bertini involution $\iota_Y$ preserves every divisor in $|{-}K_Y|$.
\item The action $\iota_Y^*\colon H^0(Y,\ol_Y(-K_Y))\to H^0(Y,\ol_Y(-K_Y))$ is $\Id$ or $-\Id$.
\item The involution $\iota_{\pr^5}$ of $\pr^5$ (resp. $\iota_Q$ of $Q$) is the identity.
\end{enumerate}
\end{rem}

Recall that we have a special surface $W_Y\subset Y$ containing $R_Y$, which is an elliptic fibration $W_Y\to \pr^1$ with fibre $F_Y$. With the same notation as in Lemma \ref{blpR}, we describe the image of $W_Y$ in $Q\subset\pr^5$.

\begin{lem}\label{imageisconic}
Every elliptic fibre $F_Y$ (resp. its strict transform $\tilde{F}_Y\subset \tilde{Y}$) is contracted by $\phi_{|{-}K_Y|}$ (resp. by $f$).
Moreover, the image of the surface $W_Y$ (resp. its strict transform $\tilde{W}_Y\subset \tilde{Y}$) is a conic $C$ in $Q\subset \pr^5$.

Furthermore, the curve $\tilde{R}_Y\coloneqq\tilde{W}_Y\cap E$ is contained in the fixed locus of $\iota_{\tilde{Y}}$, and the conic $C$ is contained in the fixed locus of $\iota_Q$.
\end{lem}
\begin{proof}
Since $-K_Y\cdot F_Y=1$, one has $\tilde{D}\cdot\tilde{F}_Y=0$, where $\tilde{D}$ is the strict transform of a general member $D\in|{-}K_Y|$. Hence $f$ contracts the elliptic fibres of $\tilde{W}_Y$. 

As $-K_Y|_{W_Y}=R_Y+2F_Y$, one has $\tilde{D}|_{\tilde{W}_Y}=(\mu^*(-K_Y)-E)|_{\tilde{W}_Y}=2\tilde{F}_Y$. Therefore, the morphism $f$ sends $\tilde{W}_Y$ to a conic in $\pr^5$.

By Lemma \ref{blpR} (i), $\iota_{\tilde{Y}}$ induces an involution on each fibre $\pr^2$ of $\mu|_E\colon E\to R_Y$. Since $W_Y$ is preserved by $\iota_Y$ by Proposition \ref{mainintroBertini}, its transform $\tilde{W}_Y\subset\tilde{Y}$ is preserved by $\iota_{\tilde{Y}}$. Therefore, the curve $\tilde{R}_Y\coloneqq\tilde{W}_Y\cap E$ (which is a section of $\mu|_E$) is invariant. Since $R_Y$ is in the fixed locus of $\iota_Y$ by Proposition \ref{mainintroBertini}, it follows that $\tilde{R}_Y$ is contained in the fixed locus of $\iota_{\tilde{Y}}$. By Lemma \ref{blpR} (ii), $f(\tilde{R}_Y)=f(\tilde{W}_Y)=C$ is contained in the fixed locus of $\iota_{Q}$.
\end{proof}

The rest of this subsection is devoted to the proofs of Theorems \ref{maininvantican} and \ref{maininvfact}. To show that $\iota_Y$ preserves every divisor in $|{-}K_Y|$, our strategy is to exclude the other remaining case by analysing the anticanonical map. 
\begin{lem}\label{directsum}
If the action $\iota_Y^*$ on  $H^0(Y,\ol_Y(-K_Y))$ is not $\pm\Id$, then
\[H^0(Y,\ol_Y(-K_Y))=V_1\oplus V_2,\]
where $V_1$ is the sub-vector space of global sections vanishing on the surface $W_Y$, and $V_2$ is uniquely determined as eigenspace corresponding to the eigenvalue $1$ or $-1$ of $\iota_Y^*$, with $\dim V_1=\dim V_2=3$. More precisely, $\iota_Y^*$ acts as $\Id$ or $-\Id$ on $V_i$ for $i=1,2$.
\end{lem}
\begin{proof}
In $\pr^4$, let $\M$ be the linear system of quintic hypersurfaces with multiplicity at least $3$ at $8$ general points. Let $\M_W$ be the sub-linear system of effective divisors in $\M$ containing the surface $W$. By ${\tt Macaulay2}$ (see Listing \ref{lst:anticanW}), one has $\mathfrak{b}(\M_W) = \mathfrak{b}(Sec)$ (see Lemma \ref{sec} for notation) for the base ideals. Moreover, by ${\tt Macaulay2}$ (see Listing \ref{lst:singQW}),  a random divisor in $\M_W$ is singular along two elliptic normal quintic curves $E_p, E_q$ through the $8$ blown up points (the two elliptic curves may coincide, in which case the divisor has multiplicity at least $3$ along this elliptic curve, and in fact the divisor is the secant variety of the elliptic curve). Moreover, there exists a unique divisor in $\M_W$ which is singular along $E_p$ and $E_q$, as $H^0(\pr^4,\mathcal{I}_W\otimes \ol_{\pr^4}(5)\otimes\mathcal{I}_{p_1,\dots,p_8}^3)\simeq H^0(\pr^1,\ol_{\pr^1}(2))\simeq\mathbb{C}^3$.

Let $V_1\subset H^0(Y,\ol_Y(-K_Y))$ be the sub-vector space of global sections vanishing on the surface $W_Y$. Let $|V_1|$ be the corresponding sub-linear system (i.e. the linear system of effective divisors in $|{-}K_Y|$ containing the surface $W_Y$). Then $\iota_Y$ preserves the family of divisors in $|V_1|$, as $\iota_Y$ preserves the surface $W_Y$ by Proposition \ref{mainintroBertini}. Since $W_Y$ is disjoint from the indeterminacy locus of the map $\eta_h\colon Y\dashrightarrow \pr^4$, and the intersection of $W_Y$ with the exceptional locus of $\eta_h$ is the union of $8$ points, we deduce that a general member in $|V_1|$ is singular along two elliptic fibres $F_{Y,1}, F_{Y,2}$ of $W_Y$, and there exists a unique divisor in $|V_1|$ which is singular along $F_{Y,1}$ and $F_{Y,2}$. Since $\iota_Y$ preserves every elliptic fibre $F_Y$ of $W_Y$ (see Proposition \ref{mainintroBertini}), we deduce that $\iota_Y$ preserves every divisor in $|V_1|$, i.e. the action of $\iota_Y^*$ on $V_1$ is $\Id$ or $-\Id$.

By Lemma \ref{restriction}(i), we have the following short exact sequence:
\[
0\to H^0(Y,\ol_{Y}(-K_Y)\otimes\mathcal{I}_{W_Y}) \to H^0(Y,\ol_Y(-K_Y))\overset{r_1}{\rightarrow} H^0(W_Y,\ol_{W_Y}(-K_Y))\to 0.
\]
Hence, $H^0(Y,\ol_Y(-K_Y))=V_1\oplus V_2$ with $V_1\simeq \Ker r_1$ and $V_2\simeq \im r_1$.

By Proposition \ref{mainintroBertini}, every elliptic fibre $F_Y$ is preserved by $\iota_Y$ and $R_Y$ is fixed by $\iota_Y$. Since $-K_Y|_{W_Y}=2F_Y+R_Y$, we deduce that $\iota_Y$ preserves every divisor in $|{-}K_Y|_{W_Y}$. Thus the action of $(\iota_Y|_{W_Y})^*$ on $H^0(W_Y,\ol_{W_Y}(-K_Y))$ is $\Id$ or $-\Id$. As $\iota_Y^*$ is not $\pm \Id$ on $H^0(Y,\ol_Y(-K_Y))$, we deduce that $V_2$ can be uniquely determined as the eigenspace corresponding to the eigenvalue $1$ or $-1$ of $\iota_Y^*$.
\end{proof}

Let $\{s_{11},s_{12},s_{13}\}$ (resp. $\{s_{21},s_{22},s_{23}\}$) be a basis of $V_1$ (resp. of $V_2$). Suppose that $\iota_Q$ is not the identity. Then for $y\in Y\backslash R_Y$, one has 
\begin{equation}\label{eqinv}
\iota_Q(\phi_{|{-}K_Y|}(y))=[s_{11}(y):s_{12}(y):s_{13}(y):-s_{21}(y):-s_{22}(y):-s_{23}(y)]
\end{equation}
by Lemma \ref{directsum}. Moreover, if $y$ is a fixed point of $\iota_Y$, then by Lemma \ref{blpR} (ii), $\phi_{|{-}K_Y|}(y)$ is fixed by $\iota_Q$. Thus by (\ref{eqinv}), one has $s_{11}(y)=s_{12}(y)=s_{13}(y)=0$ or $s_{21}(y)=s_{22}(y)=s_{23}(y)=0$, i.e. $y\in\Bs|V_1|$ or $y\in\Bs|V_2|$.

Now for $i=1,2$, let $\w{V}_i\subset H^0(\tilde{Y},\ol_{\tilde{Y}}(\mu^*(-K_Y)-E))$ be the sub-vector space of global sections which are the linear spans of $\tilde{s}_{ij}$ with $j=1,2,3$, where $\tilde{s}_{ij}$ is the strict transform of the global section $s_{ij}\in V_i$.
Hence, if a point $y\in\tilde{Y}$ is fixed by $\iota_{\tilde{Y}}$, then by repeating the argument in above paragraph, we obtain $y\in\Bs|\w{V}_1|$ or $y\in\Bs|\w{V}_2|$. To summarise, we have the following corollary.
\begin{cor}\label{fixepoints}
Suppose that $\iota_Q$ is not the identity. If a point $y\in\tilde{Y}$ is fixed by $\iota_{\tilde{Y}}$, then $y\in\Bs|\w{V}_1|$ or $y\in\Bs|\w{V}_2|$.
\end{cor}

Recall that we have the normal bundle $\ma{N}_{R_Y/Y}\cong\ol_{\pr^1}(-1)\oplus\ol_{\pr^1}^{\oplus 2}$ by Lemma \ref{normalbundle}. Hence $E=\pr(\ma{N}_{R_Y/Y}^*)\simeq\pr(\ol_{\pr^1}(1)\oplus\ol_{\pr^1}^{\oplus 2})$. Denote by $\xi$ a tautological divisor associated to $\ol_{\pr(\ma{N}_{R_Y/Y}^*)}(1)$, and $F_E$ a fibre of the projection $E=\pr(\ma{N}_{R_Y/Y}^*)\to R_Y\simeq\pr^1$. Let $l$ be an exceptional curve of $\mu$ and $\gamma$ be the curve which generates the other extremal ray $\Gamma$ of $\NE(E)$ such that $-K_E\cdot \gamma$ is the length of $\Gamma$. Then
\[
F_E\cdot l=0, \quad F_E\cdot \gamma=1,
\]
\[
\xi\cdot l=1, \quad \xi\cdot \gamma =0.
\]
Moreover, $\tilde{R}_Y\sim l+\gamma$.

With the same notation as in Lemma \ref{blpR}, we may describe the image $f(E)$ as follows.

\begin{rem}\label{imageE}
\rm{The exceptional divisor $E$ is isomorphic to the blow-up $B$ of $\pr^3$ at a line (and $B$ is embedded in $\pr^1\times\pr^3$ with bidegree $(1,1)$).

Let $\tilde{D}\subset\tilde{Y}$ be the strict transform of a general member $D\in|{-}K_Y|$. Then $(\tilde{D}|_E)^3=4$, hence $\tilde{D}|_E\sim\xi+F_E$ is very ample, with $h^0(E,\ol_E(\tilde{D}))=7$. Thus the corresponding linear system embeds $B$ in $\pr^6$ as a hyperplane section of the Segre embedding of $\pr^1\times\pr^3$ in $\pr^7$, and $B$ has degree $4$. Hence, $f|_E$ is given by the projection of $B$ from a point $x$ outside $B$ in $\pr^6$ (in fact, it is given by a sub-linear system of $|\tilde{D}|$ of dimension $5$, which is still base-point-free). 

If the point $x$ is general, then the projection is birational and the image has degree $4$ in $\pr^5$. There could be special point $x$ such that the projection has degree $2$, and the image is a $3$-dimensional quadric in $\pr^5$. In any case, the image of a fibre $F_E$ is a plane in $\pr^5$.}
\end{rem}

\begin{lem}\label{leminvQ}
Suppose that $\iota_Q$ is not the identity. Then $\iota_{\tilde{Y}}|_E$ is not the identity, and the following statements hold:
\begin{enumerate}[label=\normalfont(\roman*)]
\item The fixed locus of $\iota_{\tilde{Y}}|_E$ is the disjoint union $S_E\cup C_2$, where $S_E=\Bs|\w{V}_1|\cap E$ is the unique member in $|\xi-F_E|$ isomorphic to $\pr^1\times\pr^1$, and $C_2=\Bs|\w{V}_2|\cap E$ is a curve satisfying $C_2\sim l+\gamma$ which is mapped surjectively to $R_Y$.
\item The fixed locus of $\iota_{\pr^5}$ is two disjoint planes $\pr^2_1\cup\pr^2_2$ such that $f(S_E)=\pr^2_1$ and that $f(C_2)$ is a conic contained in $\pr^2_2$. Furthermore, $f(E)$ is a $3$-dimensional quadric in $\pr^5$.
\end{enumerate}
\end{lem}
\begin{proof}
Suppose by contradiction that $\iota_{\tilde{Y}}|_E$ is the identity. Then by Corollary \ref{fixepoints}, one has $E\subset\Bs|\w{V}_1|$ or $E\subset\Bs|\w{V}_2|$. Since $|\w{V}_i|\subset|\mu^*(-K_Y)-E|$ for $i=1,2$, this contradicts to the fact that there is no divisor in $|{-}K_Y|$ having multiplicity at least $2$ along $R_Y$ by {\tt Macaulay2} (see Listing \ref{lst:2R}).

(i) Since $\iota_{\tilde{Y}}|_{E}$ is not the identity and $\iota_Y|_{R_Y}$ is the identity, we have that $\iota_{\tilde{Y}}|_{F_E}$ is not the identity. Thus the fixed locus of $\iota_{\tilde{Y}}|_{F_E}$ is the disjoint union of a point and a line.

We first describe the base locus of $|\w{V}_1|$. Since by {\tt Macaulay2} (see Listing \ref{lst:anticanW}), one has $\mathfrak{b}(\M_W) = \mathfrak{b}(Sec)$ (see Lemma \ref{sec} for notation) for the base ideals. Thus the base locus of $|V_1|$ contains the surface $W_Y$ with multiplicity $2$ by Lemma \ref{sec}. Therefore, the base locus of $|\w{V}_1|$ contains the strict transform $\tilde{W}_Y\subset\tilde{Y}$.
Moreover, since a general member in $|V_1|$ is singular along two elliptic fibres of $W_Y$, a local computation shows that every member in $|\w{V}_1|$ contains two fibres $F_E$ above the two points on $R_Y$ where it is singular. As $\tilde{D}|_E\sim \xi+F_E$, we deduce that the unique member $S_E\in|\xi-F_E|$ is contained in the base locus of $|\w{V}_1|$.
Therefore, $\Bs|\w{V}_1|\cap E=S_E\simeq\pr^1\times\pr^1$, and $S_E\cap\tilde{W}_Y=E\cap\tilde{W}_Y=\tilde{R}_Y$.

Now we describe the base locus of $|\w{V}_2|$. Since $|\mu^*(-K_Y)-E|$ is base-point-free, $\Bs|\w{V}_2|$ is disjoint from the surface $S_E$. Let $D_2$ be a general member in $|V_2|$. Since $D_2$ does not contain the surface $W_Y$, and $D_2|_{W_Y}=R_Y+2F_Y$, we deduce that the intersection of the singular locus $\Sing D_2$ with the curve $R_Y$ contains at most one point (which is a singularity of multiplicity two). Hence a general member in $|\w{V}_2|$ contains at most one fibre $F_E$ of $E\to R_Y$.

\medskip
\noindent
\textit{Claim}. $\Bs|\w{V}_2|\cap E$ has dimension at most one.

\noindent
Suppose that there is a surface $S_2\subset \Bs|\w{V}_2|\cap E$.  Since $\tilde{D}|_E\sim \xi+F_E$, one has $S_2\in|\xi|$ or $S_2\in|\xi+F_E|$. As $\xi\cdot \tilde{R}_Y=F_E\cdot \tilde{R}_Y=1$, one has $S_2\cdot \tilde{R}_Y>0$, which contradicts the fact that $S_2$ is disjoint from $S_E$. This proves the claim.
\medskip

Note that $\Bs|\w{V}_2|\cap E$ has dimension one. This is because the fixed locus of $\iota_{\tilde{Y}}|_{F_E}$ is the disjoint union of a point and a line, and $\Bs|\w{V}_1|\cap F_E$ is a line. Thus by Corollary \ref{fixepoints}, the fixed point disjoint from the fixed line is contained in $\Bs|\w{V}_2|\cap F_E$. Therefore,  $\Bs|\w{V}_2|\cap E$ is a curve which is mapped surjectively to $R_Y$. 

Denote by $C_2$ the curve $\Bs|\w{V}_2|\cap E$. Since $C_2$ is disjoint from $S_E\in|\xi-F_E|$, one has $(\xi-F_E)\cdot C_2=0$ and thus $C_2\sim m(\gamma+l)$ with $m\geq 1$. As $(\tilde{D}|_E)^2\sim (\xi+F)^2\sim \gamma+3l$, we deduce that $C_2\sim \gamma + l$. Hence $C_2$ is an irreducible curve which is mapped surjectively to $R_Y$. 

(ii) Since $\iota_Q$ is not the identity (i.e. $\iota_{\pr^5}$ is not the identity), the fixed locus of $\iota_{\pr^5}$ is the union of two disjoint sub-linear spaces $\pr^i\cup \pr^j$ with $i+j=4$. 
Therefore, the fixed locus of $\iota_{\pr^5}$ is two disjoint planes $\pr^2_1$ and $\pr^2_2$ by the equation (\ref{eqinv}).

By Remark \ref{imageE}, $f(S_E)$ has dimension $2$. Hence, $f(S_E)\subset Q$ is one of the two planes $\pr^2_1$ and $\pr^2_2$ contained in the fixed locus of $\iota_{\pr^5}$. We may denote $f(S_E)=\pr^2_1$. Following the discussion in Remark \ref{imageE}, we now describe the map $f|_{S_E}$: the restricted linear system $|\tilde{D}|_E|_{S_E}$ embeds the surface $S_E$ in $\pr^3$ as the Segre embedding of $\pr^1\times\pr^1$ in $\pr^3$, and the image $S_E'$ has degree $2$. Hence, $f|_{S_E}$ is given by the projection of $S_E'$ from a point outside $S_E'$ in $\pr^3$. The projection has degree $2$ and the image is the plane $\pr^2_1$. Note that in Remark \ref{imageE}, the projection of $B$ from a point $x$ outside $B$ in $\pr^6$ cannot be birational, as the projection of $S_E'\subset B$ from the point $x$ in $\pr^3\subset\pr^6$ has degree $2$. Therefore, $f(E)$ is a $3$-dimensional quadric in $\pr^5$.

Since $\tilde{D}\cdot C_2=(\xi+F_E)\cdot(\gamma+l)=2$, the image $f(C_2)$ is a conic. As $f(C_2)$ is disjoint from $f(S_E)$, we deduce that $f(C_2)$ is contained in $\pr^2_2$.
\end{proof}

\begin{cor}\label{invQid}
The involution $\iota_Q$ is the identity, and thus the Bertini involution $\iota_Y$ preserves every divisor in $|{-}K_Y|$ and $f$ factors through the quotient $\tilde{Y}/\iota_{\tilde{Y}}$ via the lifted involution $\iota_{\tilde{Y}}$.
\end{cor}
\begin{proof}
Suppose by contradiction that $\iota_Q$ is not the identity. We use the notation as in Lemma \ref{leminvQ}.

Since the restricted linear system $|\tilde{D}|_E|_{S_E}$ embeds the surface $S_E$ in $\pr^3$ as the Segre embedding of $\pr^1\times\pr^1$ in $\pr^3$ with image $S_E'$ a quadric surface, the map $f|_{S_E}$ is given by the projection of $S_E'$ from a point outside $S_E'$ in $\pr^3$. The projection has degree $2$, and $f(S_E)=\pr^2_1$ by Lemma \ref{leminvQ}. Therefore, $f|_{S_E}\colon S_E\simeq \pr^1\times\pr^1\to \pr^2_1$ is a double cover branched over a non-singular conic $\Delta$ in $\pr^2_1$; moreover, the image of any line on $S_E$ is a tangent line to $\Delta$, and conversely the preimage of each tangent line on $\Delta$ is two lines on $S_E$, one from each ruling.

Let $D\in|V_1|$ be a general member and $\tilde{D}\in|\w{V}_1|$ be its strict transform. Then by Lemma \ref{leminvQ}, $\tilde{D}\cap E$ contains the surface $S_E$ and two distinct $\pr^2$ (denoted by $F_{E_1}$ and $F_{E_2}$) above the two points on $R_Y$ where $D$ is singular. Thus $f(\tilde{D})$ contains $f(F_{E_1})\eqqcolon \Pi_1$ and $f(F_{E_2})\eqqcolon\Pi_2$ which are two planes in $\pr^5$ by Remark \ref{imageE}. Moreover, $\Pi_1$ and $\Pi_2$ are distinct. This is because $F_{E_1}\cap S_E$ and $F_{E_2}\cap S_E$ are two distinct lines of a same ruling of $S_E\simeq \pr^1\times\pr^1$, and thus their images in $\pr^2_1=f(S_E)$ are two distinct tangent lines to $\Delta$ by the above discussion. Therefore, $f(\tilde{D})\cap f(E)$ contains three distinct planes $\pr^2_1,\Pi_1,\Pi_2$. This contradicts to the fact that $f(\tilde{D})$ is a hyperplane in $\pr^5$ and $f(E)$ is a $3$-dimensional quadric in $\pr^5$ (so that their intersection is a surface of degree $2$ in $\pr^5$).

Therefore, $\iota_Q$ is the identity. By Lemma \ref{blpR} (ii), one has $f=f\circ\iota_{\tilde{Y}}$. Thus $f$ factors through the quotient $\tilde{Y}/\iota_{\tilde{Y}}$.
\end{proof}

\begin{cor}\label{invE}
The morphism $f|_E\colon E\to f(E)$ is birational, and $f(E)$ has degree $4$ in $\pr^5$. Moreover, the restricted involution $\iota_{\tilde{Y}}|_E$ is the identity.
\end{cor}
\begin{proof}
By Remark \ref{imageE}, either $f|_E$ has degree $2$ and the image is a $3$-dimensional quadric in $\pr^5$, or $f|_E$ is finite birational and the image has degree $4$ in $\pr^5$. We will show that the first case cannot happen.

Suppose that $f|_E$ has degree $2$ and $f(E)$ is a $3$-dimensional quadric in $\pr^5$. We will follow the same argument as in the proof of Corollary \ref{invQid}. By Lemma \ref{leminvQ} (with the same notation), for a general member $D\in|V_1|$, its strict transform  $\tilde{D}\in|\w{V}_1|$ contains the surface $S_E\subset E$. Moreover, $\tilde{D}$ contains the two distinct fibres (denoted by $F_{E_1}$ and $F_{E_2}$) of $\mu|_E\colon E\to R_Y$ above the two points on $R_Y$ where $D$ is singular. Hence, $f(\tilde{D})\cap f(E)$ contains the surface $f(S_E)$ and the two planes $f(F_{E_1}),f(F_{E_2})$ (which may coincide).
\begin{enumerate}[label=\normalfont(\alph*)]
\item If $f|_{S_E}$ has degree $2$, then the same argument as in the proof of Corollary \ref{invQid} shows that $f(S_E),f(F_{E_1}),f(F_{E_2})$ are $3$ distinct planes.
\item If $f|_{S_E}$ has degree $1$, then $f(S_E)$ is either a non-normal surface or isomorphic to $S_E\simeq \pr^1\times\pr^1$. Thus $f(S_E)$ has degree at least $2$ in $\pr^5$.
\end{enumerate}
This contradicts to the fact that $f(\tilde{D})\cap f(E)$ is a surface of degree $2$ in $\pr^5$. Therefore, $f|_E$ is finite birational and $f(E)$ has degree $4$ in $\pr^5$.

Now suppose that $\iota_{\tilde{Y}}|_E$ is not the identity. Since $f|_E=f|_E\circ\iota_{\tilde{Y}}|_E$ by Corollary \ref{invQid}, we deduce that $f|_E$ has degree $2$, which leads to a contradiction.
\end{proof}

\begin{proof}[Proof of Theorem \ref{maininvantican}]
Follows from Corollary \ref{invQid}.
\end{proof}

\begin{proof}[Proof of Theorem \ref{maininvfact}]
Follows from Lemmas \ref{blpR}, \ref{imageisconic}, and Corollaries \ref{invQid}, \ref{invE}.
\end{proof}

\begin{rem}\label{fixedlocustildeY}
The fixed locus of $\iota_{\tilde{Y}}$ is $E\cup Res$, where $Res$ has dimension at most $2$ and its intersection with every $\tilde{P}_{\ell}$ is non-empty and zero-dimensional, where $\tilde{P}_{\ell}\subset\tilde{Y}$ is the strict transform of $P_\ell$ (see notation in Proposition \ref{smallloci}).
\end{rem}
\begin{proof}
Let $P_\ell\simeq \pr^2$ be the exceptional locus of a small extremal contraction of $Y$. Then $\iota_Y(P_\ell)=P_{\iota_S^*(\ell)}$ is also the exceptional locus of some small extremal contraction of $Y$ and $P_\ell$ intersects $\iota_Y(P_\ell)$ transversally at $3$ points by \cite[Rem.~2.15 (c), Lem.~6.4]{MR3942925}. Therefore, the intersection of $P_\ell$ with the fixed locus of $\iota_Y$ is non-empty and zero-dimensional. 

As $R_Y\subset W_Y$ is disjoint from $P_\ell$ by Lemma \ref{Pldisj}, we deduce that $E$ is disjoint from $\tilde{P}_{\ell}$. Therefore, $Res\cap\tilde{P}_\ell$ is non-empty and zero-dimensional. As every non-zero effective divisor in $Y$ must have positive intersection with some extremal ray of $\NE(Y)$, we deduce that $Res$ has dimension at most $2$.
\end{proof}

\newpage
\begin{appendix}
\section{Computations by {\tt Macaulay2}}\label{appen}
\begin{lstlisting}
restart
k = ZZ/67
\end{lstlisting}
We set up the projective space $\pr^4$:
\begin{lstlisting}
R = k[x_0..x_4]
\end{lstlisting}
We choose $8$ points in $\pr^4$:
\begin{lstlisting}
I_0 = ideal(x_1,x_2,x_3,x_4)
I_1 = ideal(x_0,x_2,x_3,x_4)
I_2 = ideal(x_1,x_0,x_3,x_4)
I_3 = ideal(x_1,x_2,x_0,x_4)
I_4 = ideal(x_1,x_2,x_3,x_0)
I_5 = ideal(x_1-x_2,x_2-x_3,x_3-x_4,x_0-x_4)
I_6 = ideal(x_0-3*x_1,x_1-7*x_2,x_2-11*x_3,x_3-13*x_4)
I_7 = ideal(x_0-17*x_1,x_1-23*x_2,x_2-29*x_3,x_3-31*x_4)
\end{lstlisting}
We compute the ideal \texttt{II} defined by the $6$ quintics through the $8$ points with multiplicity at least $3$:
\begin{lstlisting}
J = I_0; for j from 1 to 7 do J = intersect(J,I_j);
H = saturate J^3;
G = gens (H);
betti G
G1 = submatrix (G,{0..5});
II = ideal(G1);
III =  sheaf module II;
HH^0(III(5))
\end{lstlisting}
We check that \texttt{II} is the intersection of the ideal of the $28$ lines, the ideal of the $8$ quartics and the ideal \texttt{I5} of a smooth rational quintic curve:
\begin{lstlisting}[label=lst:base,caption=Base scheme,captionpos=b]
LL = ideal(1_R); 
for i from 0 to 7 do for j from 0 to i-1 do
LL = intersect(LL,ideal submatrix(gens intersect (I_i,I_j),{0..2}));
isSubset(II,LL)
RN = ideal(1_R);
for i from 0 to 7 do 
RN = intersect(RN,minors(2,submatrix((res (J:I_i)).dd_4,{3..6},{0..1})));
isSubset(II,RN)
I5 = ((II:LL):RN);
degree I5, genus I5, ideal singularLocus variety I5
II == intersect(intersect(LL,RN),I5)
\end{lstlisting}
We compute the normal bundle of the smooth rational quintic curve:
\begin{lstlisting}[label=lst:normal,caption=Normal bundle,captionpos=b]
RI5 = R/I5
N5 = (module I5)**RI5
PI5 = Proj RI5
SN5 = sheaf N5
HH^0(SN5)
HH^0(sheaf dual N5)
KI5 = Ext^3(R^1/I5,R^{-5})**RI5
HH^0(SN5**OO_PI5(1)**(sheaf dual KI5))
\end{lstlisting}
We choose three points on the smooth rational quintic curve:
\begin{lstlisting}
P1 = ideal(x_3-14*x_4,x_2-x_4,x_1+x_4,x_0-12*x_4)
P2 = ideal(x_3+17*x_4,x_2-22*x_4,x_1+20*x_4,x_0+2*x_4)
P3 = ideal(x_3-26*x_4,x_2+27*x_4,x_1-30*x_4,x_0+21*x_4)
\end{lstlisting}
We compute the quintic with multiplicity $3$ at the $8$ points and the point \texttt{P1} (resp. \texttt{P2} and resp. \texttt{P3}):
\begin{lstlisting}[label=lst:secant,caption=Three secant varieties,captionpos=b]
J13 = intersect(J^3,P1^3);
H13 = saturate J13;
G13 = gens(H13);
betti G13
GP1 = submatrix(G13,{0});
Q1 = ideal(GP1);

J23 = intersect(J^3,P2^3);
H23 = saturate J23;
G23 = gens(H23);
betti G23
GP2 = submatrix(G23,{0});
Q2 = ideal(GP2);

J33 = intersect(J^3,P3^3);
H33 = saturate J33;
G33 = gens(H33);
betti G33
GP3 = submatrix(G33,{0});
Q3 = ideal(GP3);
\end{lstlisting}
We compute the elliptic normal quintic curve along which \texttt{Q1} is singular:
\begin{lstlisting}
SingQ1 = ideal singularLocus variety Q1;
dim SingQ1, degree SingQ1
SSingQ1 = ideal singularLocus variety SingQ1;
dim SSingQ1, degree SSingQ1
E1 = (SingQ1:SSingQ1);
dim E1, degree E1, genus E1
ideal singularLocus variety E1
\end{lstlisting}
We compute the intersection of the three quintics and obtain the cubic scroll \texttt{W}:
\begin{lstlisting}[label=lst:W,caption=Scheme-theoretic intersection of secant varieties,captionpos=b]
SS3 = Q1 + Q2 + Q3;
SS = (SS3:II);
dim SS, degree SS
W = ideal singularLocus variety SS;
dim W, degree W, genus W, ideal singularLocus variety W
W == (SS:W)
\end{lstlisting}
We compute the quintics through the $8$ points with multiplicity at least $3$ containing the surface \texttt{W}:
\begin{lstlisting}[label=lst:anticanW,caption=Quintics containg $W$,captionpos=b]
JW = intersect(J^3,W);
HW = saturate JW;
GW = gens(HW);
betti GW
GW1 = submatrix(GW,{0..2});
IIW = ideal(GW1);
IIIW = sheaf module IIW;
HH^0(IIIW(5))
IIW == SS3
\end{lstlisting}
We look at the singular locus of a quintic hypersurface through the $8$ points with multiplicity at least $3$ containing the surface \texttt{W}:
\begin{lstlisting}[label=lst:singQW,caption=Singular locus of a quintic containg $W$,captionpos=b]
QW =  ideal(11*GW1_(0,0)+7*GW1_(0,1)+19*GW1_(0,2));
SingQW = ideal singularLocus variety QW;
dim SingQW, degree SingQW, genus SingQW
SingW1 = (SingQW:J);
SingW2 = (SingW1:J);
SingW3 = (SingW2:J);
LSingQW = decompose SingW3;
EE = ideal(LSingQW);
degree EE, genus EE, dim EE
SingEE = ideal singularLocus variety EE;
dim SingEE, degree SingEE
degree (E1 + EE)
degree (I5 + EE)
\end{lstlisting}
We compute the hypersurfaces of degree $10$ through the $8$ points with multiplicity at least $6$:
\begin{lstlisting}
JI5 = intersect(I_5^6,intersect(I_6^6,I_7^6));
HI5 = saturate JI5;
GI5 = gens HI5;
betti GI5
GI15 = submatrix(GI5,{0..105});
II5 = ideal(GI15);
JI4 = intersect(I_4^6,II5);
HI4 = saturate JI4;
GI4 = gens HI4;
betti GI4
GI14 = submatrix(GI4,{0..132});
II4 = ideal(GI14);
JI3 = intersect(I_3^6,II4);
HI3 = saturate JI3;
GI3 = gens HI3;
betti GI3
GI13 = submatrix(GI3,{0..154});
II3 = ideal(GI13);
JI2 = intersect(I_2^6,II3);
HI2 = saturate JI2;
GI2 = gens HI2;
betti GI2
GI12 = submatrix(GI2,{0..123});
II2 = ideal(GI12);
JI1 = intersect(I_1^6,II2);
HI1 = saturate JI1;
GI1 = gens HI1;
betti GI1
GI11=submatrix(GI1,{0..136});
II1 = ideal(GI11);
JI0 = intersect(I_0^6,II1);
HI0 = saturate JI0;
GI0 = gens HI0;
betti GI0
GG = submatrix(GI0,{0..28});
IGG = ideal(GG);
\end{lstlisting}
We compute the hypersurfaces of degree $10$ through the $8$ points with multiplicity at least $6$ containing the surface \texttt{W}:
\begin{lstlisting}[label=lst:bianticanW,caption=Decics containing $W$,captionpos=b]
JW2 = intersect(W,IGG);
HW2 = saturate JW2;
GW2 = gens HW2;
betti GW2
GGW = submatrix(GW2,{0..20});
IW2 = ideal(GGW);
IIW2 = sheaf module IW2;
HH^0(IIW2(10))
\end{lstlisting}
We compute the image of the elliptic normal quintic \texttt{E1} via the map defined by the linear system of hypersurfaces of degree $10$ through the $8$ points with multiplicity at least $6$:
\begin{lstlisting}
S2 = k[u_0..u_28];
ImE1 = ker map(R/E1,S2,GG);
dim ImE1, degree ImE1
\end{lstlisting}
We compute the image of the rational quintic curve \texttt{I5} via the map defined by the linear system of hypersurfaces of degree $10$ through the $8$ points with multiplicity at least $6$:
\begin{lstlisting}
ImI5 = ker map(R/I5,S2,GG);
dim ImI5, degree ImI5
\end{lstlisting}
We compute the image of the surface \texttt{W} via the map defined by the linear system of hypersurfaces of degree $10$ through the $8$ points with multiplicity at least $6$:
\begin{lstlisting}[label=lst:image,caption=Some images by the bianticanonical map,captionpos=b]
ImW = ker map(R/W, S2,GG);
dim ImW, degree ImW
\end{lstlisting}
We compute the hypersurfaces of degree $10$ with multiplicity at least $7$ at the point \texttt{I\_0} and multiplicity at least $6$ at the other $7$ points:
\begin{lstlisting}
JI00 = intersect(I_0^7,II1);
HI00 = saturate JI00;
GI00 = gens HI00;
betti GI00
GG0 = submatrix(GI00,{0});
IGG0 = ideal(GG0);
\end{lstlisting}
And we obtain a unique such hypersurface of degree $10$; now we check if this hypersurface contains the surface \texttt{W}:
\begin{lstlisting}[label=lst:special,caption=Special member in the bianticanonical system,captionpos=b]
JW0 = intersect(W,IGG0);
HW0 = saturate JW0;
GW0 = gens HW0;
betti GW0
\end{lstlisting}
We compute the image of $\pr^4$ via the map defined by the linear system of quintic hypersurfaces through the $8$ points with multiplicity at least $3$.
\begin{lstlisting}[label=lst:quadric, caption=Image by the anticanonical map,captionpos=b]
JJ = minors(2,random(R^{4:0},R^{-2,-3}));
degree JJ
genus JJ
betti res JJ == betti res ideal(G1)
S = k[y_0..y_5] ; g = map(R,S,gens JJ);
K = ker g;
dim K
degree K
singularLocus variety K
\end{lstlisting}
We check that there is no quintic hypersurfaces through the $8$ points with multipilicity at least $3$ and having multiplicity at least $2$ along the smooth rational quintic curve \texttt{I5}:
\begin{lstlisting}[label=lst:2R, caption=Quintics having multiplicity 2 along the smooth rational quintic base curve,captionpos=b]
JRR = intersect(J^3,I5^2);
HRR = saturate JRR;
GRR = gens HRR;
betti GRR
\end{lstlisting}
\end{appendix}
\newpage
\bibliographystyle{alpha}
  \bibliography{bibliography.bib}
  
  \Affilfont{\small{Z\textsc{hixin} X\textsc{ie}, F\textsc{achrichtung} M\textsc{athematik}, C\textsc{ampus}, G\textsc{eb\"{a}ude} E2.4, U\textsc{niversit\"{a}t} \textsc{des} S\textsc{aarlandes}, 66123 S\textsc{aarbr\"{u}cken}, G\textsc{ermany} }}

\textit{Email address:}
\href{mailto:xie@math.uni-sb.de}{xie@math.uni-sb.de}
\end{document}